\documentclass[final,onefignum,onetabnum,letterpaper]{siamonline250211}
\usepackage[scale=0.8]{geometry} 

\usepackage{soul}
\usepackage{lipsum}
\usepackage{amsfonts}
\usepackage{epstopdf}
\ifpdf
  \DeclareGraphicsExtensions{.eps,.pdf,.png,.jpg}
\else
  \DeclareGraphicsExtensions{.eps}
\fi

\usepackage{datetime}
\newdateformat{monthyeardate}{%
	\monthname[\THEMONTH] \THEDAY, \THEYEAR}

\usepackage{academicons}
\usepackage{xcolor}
\renewcommand{\orcid}[1]{\href{https://orcid.org/#1}{\textcolor[HTML]{A6CE39}{orcid.org/#1}}}

\usepackage{amsmath}
\allowdisplaybreaks
\usepackage{amssymb}
\usepackage{commath}
\usepackage{mathtools}
\usepackage{bbm}
\usepackage{bm}

\usepackage{color}
\usepackage{graphicx}
\usepackage[small]{caption}
\usepackage{subcaption}

\usepackage{relsize}
\usepackage{adjustbox}
\usepackage{algorithm}
\usepackage[noend]{algpseudocode}
\usepackage{booktabs}
\usepackage{tikz}
\usepackage{cancel}
\usepackage{pifont}
\usepackage{multirow} 
\usetikzlibrary{bayesnet} 
\usetikzlibrary{decorations.pathreplacing} 
\usetikzlibrary{arrows.meta} 

\usepackage[nocompress]{cite}

\usepackage{enumitem}
\setlist[enumerate]{leftmargin=.5in}
\setlist[itemize]{leftmargin=.5in}

\newsiamremark{remark}{Remark}
\newsiamremark{example}{Example}
\newsiamremark{hypothesis}{Hypothesis}
\newsiamremark{assumption}{Assumption}
\newsiamthm{claim}{Claim}

\AddToHook{env/remark/begin}{\crefalias{theorem}{remark}}
\AddToHook{env/example/begin}{\crefalias{theorem}{example}}
\AddToHook{env/hypothesis/begin}{\crefalias{theorem}{hypothesis}}
\AddToHook{env/assumption/begin}{\crefalias{theorem}{assumption}}

\crefname{remark}{Remark}{Remarks}
\crefname{example}{Example}{Examples}
\crefname{hypothesis}{Hypothesis}{Hypotheses}
\crefname{assumption}{Assumption}{Assumptions}

\headers{Edge-preserving BNNs via hierarchical priors}{J.\ Habscheid and J.\ Glaubitz}

\title{Edge-preserving Bayesian neural networks via hierarchical sparsity-promoting priors\thanks{\monthyeardate\today 
\corresponding{Jan Glaubitz} 
}}
\author{ 
Jan Habscheid\thanks{
Division of Applied Mathematics, Link\"oping University, Sweden
(\email{jan.habscheid@liu.se}, \orcid{0009-0005-0040-1562}; 
\email{jan.glaubitz@liu.se}, \orcid{0000-0002-3434-5563})
} 
\and
Jan Glaubitz\footnotemark[2]
}

\usepackage{amsopn}

\usepackage{comment}
\definecolor{blue1}{RGB}{0,76,153}
\definecolor{green1}{RGB}{1,109,57}
\definecolor{red1}{RGB}{153,0,0}
\definecolor{orange1}{RGB}{255,127,80}

\usepackage[normalem]{ulem}

\usepackage[normalem]{ulem}

\newcommand{\MAP}{\mathrm{MAP}}

\DeclareMathOperator{\diag}{diag}

\DeclareMathOperator*{\argmin}{arg\,min}

\newcommand{\N}{\mathbb{N}}

\newcommand{\R}{\mathbb{R}}

\usepackage{lineno}
\usepackage{comment}

\ifpdf
\hypersetup{
	pdftitle={GSBL-BNNs},
 	pdfauthor={Habscheid and Glaubitz}
}

\fi

\begin{document}

\maketitle

\begin{abstract}
	Recovering piecewise-smooth functions with jump discontinuities from noisy, indirect, and incomplete data is a central task in computational inverse problems.
Bayesian neural networks (BNNs) with heavy-tailed weight priors have emerged as effective edge-preserving priors for such problems.
However, posterior inference for these models remains computationally challenging, and existing approaches settle for point estimates obtained by directly optimizing heavy-tailed log-densities.
We propose GSBL-BNNs, which model the neural network parameters hierarchically using sparsity-promoting priors from generalized sparse Bayesian learning (GSBL). 
The hierarchical formulation recovers several popular marginal priors, including Laplace and Student's $t$ (including Cauchy).
A key modeling ingredient, novel among edge-preserving BNNs and native to the hierarchical formulation, couples each neuron's bias to its weights through a shared hyperparameter, so that the induced jump discontinuities are free to occur anywhere in the domain while the conditionally Gaussian structure underlying efficient inference is preserved.
Furthermore, we develop a block-coordinate descent (BCD) algorithm for maximum a posteriori (MAP) estimation. 
This BCD approach alternates between standard network training with iteratively reweighted $\ell^2$-regularization and exact, closed-form hyperparameter updates for gamma and inverse gamma hyperpriors. 
Building on the BCD algorithm, we perform approximate Bayesian inference and uncertainty quantification (UQ) using deep ensembles. 
Computational experiments---including signal denoising, deblurring, and inverting a two-dimensional hyperbolic conservation law--- show that GSBL-BNNs deliver accurate, edge-preserving reconstructions with built-in UQ.
Compared with Gaussian- and Cauchy-BNNs, they yield more accurate predictive means and avoid overconfident uncertainty estimates.
\end{abstract}



\begin{Code}
    \url{https://doi.org/10.5281/zenodo.23184383}
\end{Code}

\begin{DOI}
    Not yet assigned
\end{DOI}

\section{Introduction}
\label{sec:intro}

Many applications in science and engineering require recovering an unknown function $u$ from noisy, indirect, and incomplete observational data $\mathbf{y}$.
Prominent examples include signal and image recovery and estimating states of hyperbolic conservation laws, where the sought-after function is often piecewise smooth with jump discontinuities.
The Bayesian approach \cite{stuart2010inverse,calvetti2023bayesian} addresses the ill-posedness of such inverse problems by combining the observational data with a prior model, resulting in a posterior distribution that also quantifies uncertainty.
Specifically, given a likelihood $f(u; \mathbf{y})$ that quantifies how well a candidate $u$ explains the data $\mathbf{y}$ and a prior density $\pi^0(u)$ that reflects our assumptions on $u$ before seeing the data, Bayes' formula combines the two into the posterior density $\pi^{\mathbf{y}}$:
\begin{equation}\label{eq:Bayes}
	\pi^{\mathbf{y}}(u)
	\propto f(u; \mathbf{y})\, \pi^0(u).
\end{equation}

Most existing \emph{edge-preserving} priors---priors whose realizations admit sharp features---act on a fixed discretization and promote sparsity of a linear transform, such as increments or gradients.
One line of work places heavy-tailed distributions---including Laplace, Cauchy, and Besov models---directly on the transformed unknown \cite{babacan2009bayesian,markkanen2019cauchy,suuronen2022cauchy,uribe2022hybrid,lan2025spatiotemporal} or employs global-local shrinkage priors such as the horseshoe \cite{dong2023inducing,uribe2023horseshoe}.
Another line of work, which the present paper builds on, instead models the transformed unknown hierarchically, combining conditionally Gaussian priors with gamma-type hyperpriors  \cite{tipping2001sparse,calvetti2020sparse,glaubitz2023generalized,xiao2023sequential,glaubitz2024leveraging,lindbloom2025efficient,glaubitz2025efficient}.
The hierarchical structure of these models admits particularly efficient inference algorithms, and generalized-gamma-based models have been found to accurately describe natural image statistics \cite{marks2025generalized}.

An emerging alternative is to parametrize the unknown function by a neural network and place the prior on the network parameters, resulting in a Bayesian neural network (BNN) \cite{izmailov2021bayesian,tran2022all,arbel2026primer,goan2026bayesian}.
BNN priors are mesh-free and expressive.
Crucially, the choice of the weight prior determines the character of the induced function-space prior: Gaussian weight priors produce smooth, Gaussian-process-like realizations that smear out edges \cite{neal1996bayesian}, while \emph{heavy-tailed} weight priors produce BNN realizations with jump-like features \cite{li2022bayesian,lee2023deep,suuronen2023bayesian,soto2025wide}.

At the same time, posterior inference for BNNs with heavy-tailed priors remains a serious computational challenge \cite{izmailov2021bayesian}, as the corresponding posteriors are high-dimensional, multimodal, and heavy-tailed.
Existing works therefore settle for maximum a posteriori (MAP) estimation by directly optimizing the heavy-tailed log-densities.
This has two drawbacks:
First, heavy-tailed log-densities have vanishing gradients in their tails, yielding flat optimization landscapes, and some (e.g., $\alpha$-stable) densities must additionally be approximated numerically \cite{suuronen2023bayesian}.
Second, a single MAP estimate provides no uncertainty quantification (UQ).

\subsection*{Our contribution}

\begin{figure}[tb]
	\centering
	\includegraphics[width=\textwidth]{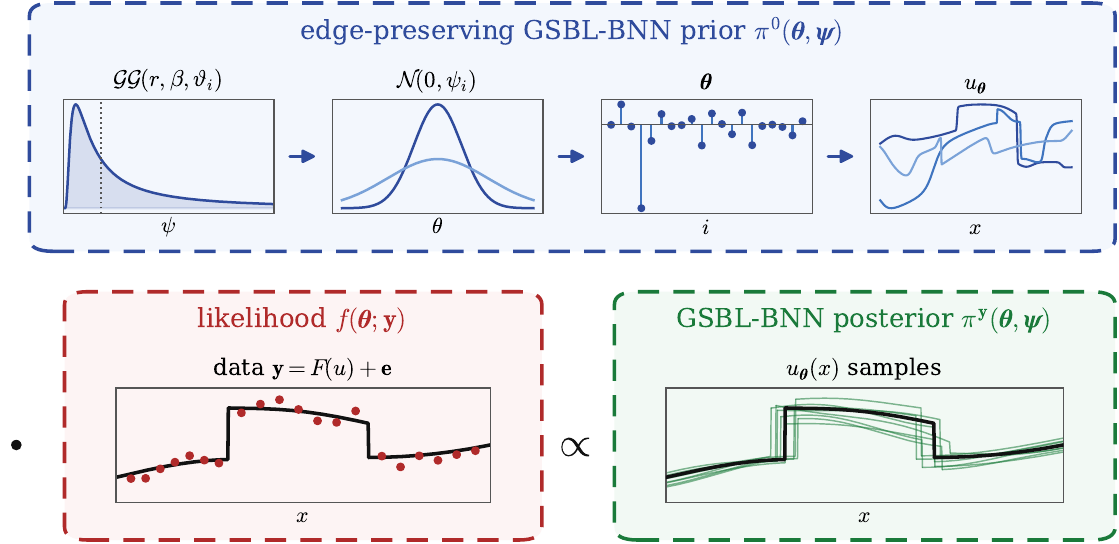}
	\caption{
		Schematic of the proposed GSBL-BNN approach.
		\emph{Top:} The \textcolor{blue1}{edge-preserving GSBL-BNN prior} (\textcolor{blue1}{blue}) is built by drawing each variance $\psi_i$ from a heavy-tailed hyperprior $\mathcal{GG}(r, \beta, \vartheta_i)$, drawing the conditionally Gaussian neural parameter $\theta_i \mid \psi_i \sim \mathcal{N}(0, \psi_i)$, and evaluating the network $u_{\boldsymbol{\theta}}$.
		\emph{Bottom:} Combining the prior with the \textcolor{red1}{likelihhood} (\textcolor{red1}{red}) via Bayes' theorem yields the \textcolor{green1}{GSBL-BNN posterior} (\textcolor{green1}{green}), which provides edge-preserving reconstructions.
	}
	\label{fig:pipeline}
\end{figure}

We address both drawbacks by modeling the neural network parameters \emph{hierarchically}, using sparsity-promoting priors from generalized sparse Bayesian learning (GSBL) \cite{calvetti2020sparse,glaubitz2023generalized,glaubitz2025efficient}. 
These combine zero-mean conditionally Gaussian priors for the neural parameters with generalized gamma hyperpriors for their variances.
We call the resulting model class \emph{GSBL-BNNs} (illustrated in \cref{fig:pipeline}). 
Our main contributions are as follows:
\begin{itemize}
	\item 
	\emph{Model} (\Cref{sec:BNNs}):
	We introduce GSBL-BNNs, which recover an exact hierarchical reformulation of heavy-tailed BNN priors.
	As a key modeling ingredient, new among edge-preserving BNNs, we couple each neuron's bias to its weights through a shared hyperparameter. 
	This frees the induced jumps to occur anywhere in the domain while preserving the conditionally Gaussian structure our inference exploits; existing edge-preserving BNNs \cite{li2022bayesian} instead draw biases independently, concentrating the jumps near the center of the domain.

	\item 
	\emph{MAP estimation} (\Cref{sec:MAP}):
	We propose a block-coordinate descent (BCD) algorithm for joint MAP estimation that alternates between standard network training with iteratively reweighted $\ell^2$-regularization and exact hyperparameter updates, available in closed form for gamma and inverse gamma hyperpriors. 
	Our numerical experiments indicate that the BCD approach makes training GSBL-BNNs more robust compared to joint optimization of the neural network parameters and hyperparameters. 

	\item 
	\emph{Uncertainty quantification} (\Cref{sec:deepEnsembles}):
	We perform approximate Bayesian inference using deep ensembles.
	We validate the resulting uncertainty estimates by their frequentist coverage.  

	\item \emph{Computational experiments} (\Cref{sec:numerics}):
	We demonstrate GSBL-BNNs for signal denoising, deblurring, and inverting a two-dimensional hyperbolic conservation law, including comparisons with existing Gaussian- and Cauchy-BNNs.
	
\end{itemize}
Our proposed GSBL-BNN approach avoids density approximations and replaces flat heavy-tailed optimization landscapes with alternating smooth subproblems. 
Notably, we focus on MAP estimation and deep ensembles. 
Other inference strategies for GSBL-BNNs, including Markov chain Monte Carlo sampling and variational inference, are the subject of future work.

\subsection*{Outline}

\Cref{sec:BNNs} introduces GSBL-BNNs.
\Cref{sec:MAP} develops the BCD approach for MAP estimation.
\Cref{sec:deepEnsembles} discusses deep ensembles for approximate Bayesian inference and UQ.
\Cref{sec:numerics} presents computational experiments, and \Cref{sec:summary} offers concluding thoughts. 
\section{GSBL-BNNs}
\label{sec:BNNs}

We now introduce the proposed GSBL-BNN model.

\subsection{Neural-network parametrization and likelihood}
\label{subsec:setup}

Consider recovering an unknown function $u\colon \Omega \to \R$ on a bounded domain $\Omega \subset \R^d$ from observational data:
\begin{equation}\label{eq:data_model}
	\mathbf{y} = F(u) + \mathbf{e},
\end{equation}
where $\mathbf{y} \in \R^{M}$ is the given observational data, $F$ is a known (possibly nonlinear) forward operator, and $\mathbf{e} \in \R^{M}$ is observational noise.
We are interested in piecewise-smooth functions $u$ with jump discontinuities. 
These arise, for instance, in signal and image recovery \cite{kaipio2006statistical,calvetti2007gaussian,suuronen2022cauchy,glaubitz2023generalized} and as solutions of hyperbolic conservation laws \cite{lax1973hyperbolic,dafermos2005hyberbolic}, where locating the jumps is also central to numerical shock capturing \cite{krivodonova2004shock,ray2018artificial,glaubitz2019high}.

\begin{figure}[tb]
	\centering
	\resizebox{0.8\textwidth}{!}{
\begin{tikzpicture}[
	neuron/.style={circle, draw, minimum size=0.95cm, inner sep=1pt, font=\small},
	inputnode/.style={neuron, draw=green!55!black, fill=green!12},
	hiddenA/.style={neuron, draw=blue!55!black, fill=blue!10},
	hiddenB/.style={neuron, draw=orange!75!black, fill=yellow!22},
	outputnode/.style={neuron, draw=red!65!black, fill=red!10},
	conn/.style={->, semithick, black!65},
	headerlab/.style={font=\footnotesize, align=center},
	bracelab/.style={font=\small}
]
	\def\xin{0} \def\xhA{2.5} \def\xhB{5.0} \def\xout{7.5}

	\node[headerlab] at (\xin, 2.9) {Input layer};
	\node[headerlab] at (\xhA, 2.9) {Hidden layer 1\\ (width $n_1$)};
	\node[headerlab] at (\xhB, 2.9) {Hidden layer 2\\ (width $n_2$)};
	\node[headerlab] at (\xout, 2.9) {Output layer};

	\node[inputnode] (x) at (\xin, 0) {$\boldsymbol{x}$};
	\draw[conn, green!55!black] (-1.5, 0) node[left, font=\small, green!55!black] {Input} -- (x);

	\node[hiddenA] (a1) at (\xhA, 1.8) {$h_1^{(1)}$};
	\node[hiddenA] (a2) at (\xhA, 0.6) {$h_2^{(1)}$};
	\node          (adots) at (\xhA, -0.55) {$\vdots$};
	\node[hiddenA] (a3) at (\xhA, -1.8) {$h_{n_1}^{(1)}$};

	\node[hiddenB] (b1) at (\xhB, 1.8) {$h_1^{(2)}$};
	\node[hiddenB] (b2) at (\xhB, 0.6) {$h_2^{(2)}$};
	\node          (bdots) at (\xhB, -0.55) {$\vdots$};
	\node[hiddenB] (b3) at (\xhB, -1.8) {$h_{n_2}^{(2)}$};

	\node[outputnode] (out) at (\xout, 0) {};
	\draw[conn, red!65!black] (out) -- (9, 0) node[right, font=\small, red!65!black] {$u_{\boldsymbol{\theta}}(\boldsymbol{x})$};

	\foreach \t in {a1, a2, a3} { \draw[conn] (x) -- (\t); }
	\foreach \s in {a1, a2, a3} { \foreach \t in {b1, b2, b3} { \draw[conn] (\s) -- (\t); } }
	\foreach \s in {b1, b2, b3} { \draw[conn] (\s) -- (out); }

	\draw[decorate, decoration={brace, mirror, amplitude=5pt}]
		(\xin + 0.1, -2.6) -- (\xhA - 0.1, -2.6)
		node[bracelab, midway, below=6pt] {$W^{(1)}, \boldsymbol{b}^{(1)}$};
	\draw[decorate, decoration={brace, mirror, amplitude=5pt}]
		(\xhA + 0.1, -2.6) -- (\xhB - 0.1, -2.6)
		node[bracelab, midway, below=6pt] {$W^{(2)}, \boldsymbol{b}^{(2)}$};
	\draw[decorate, decoration={brace, mirror, amplitude=5pt}]
		(\xhB + 0.1, -2.6) -- (\xout - 0.1, -2.6)
		node[bracelab, midway, below=6pt] {$W^{(3)}, b^{(3)}$};
\end{tikzpicture}
	}
	\caption{
		Fully-connected feedforward neural network \cref{eq:NN}, mapping $\boldsymbol{x} \in \Omega$ to $u_{\boldsymbol{\theta}}(\boldsymbol{x}) \approx u$. 
		The hidden unit $h_i^{(\ell)}$ is the $i$th component of the layer activation $\boldsymbol{h}^{(\ell)} = \phi\bigl( W^{(\ell)} \boldsymbol{h}^{(\ell-1)} + \boldsymbol{b}^{(\ell)} \bigr)$, with $\boldsymbol{h}^{(0)} = \boldsymbol{x}$.
	}
	\label{fig:NN_schematic}
\end{figure}
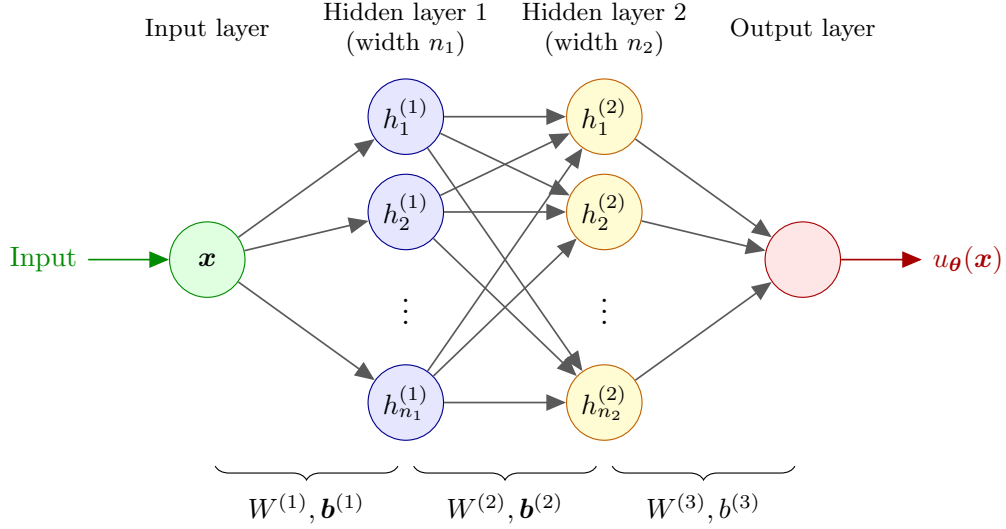

Following \cite{li2022bayesian}, we parametrize the unknown function $u$ by a fully-connected feedforward neural network $u_{\boldsymbol{\theta}}\colon \Omega \to \R$ with $L$ layers:
\begin{equation}\label{eq:NN}
	u_{\boldsymbol{\theta}}(\boldsymbol{x})
		= W^{(L)} \phi\Bigl( \cdots \phi\bigl( W^{(2)} \overbrace{\phi\bigl( \underbrace{W^{(1)} \boldsymbol{x} + \boldsymbol{b}^{(1)}}_{\text{pre-activation } \boldsymbol{z}^{(1)}} \bigr)}^{\text{activation } \boldsymbol{h}^{(1)}} + \boldsymbol{b}^{(2)} \bigr) \cdots \Bigr) + b^{(L)},
\end{equation}
where $W^{(\ell)} \in \R^{n_{\ell} \times n_{\ell-1}}$ and $\boldsymbol{b}^{(\ell)} \in \R^{n_{\ell}}$ are the weight matrices and bias vectors, $n_\ell \in \N$ is the width of the $\ell$th layer with $n_0 = d$ and $n_L = 1$, and the activation function $\phi$ is applied componentwise.
We collect all weights and biases—the \emph{neural parameters}—in a single vector $\boldsymbol{\theta} \in \R^K$.
\Cref{fig:NN_schematic} illustrates the neural network structure.
Throughout, we use the bounded and smooth activation $\phi = \tanh$, consistent with \cite{li2022bayesian,soto2025wide}; also see \cref{rem:tanh} below. 
Notably, the computational goal is to infer the neural parameters $\boldsymbol{\theta}$. 

For simplicity, we model the noise in \cref{eq:data_model} as i.i.d.\ zero-mean Gaussian, $\mathbf{e} \sim \mathcal{N}(\mathbf{0}, \sigma^2 I_{M})$, with known variance $\sigma^2$; see \cref{rem:general_noise} for generalizations.
Substituting $u \approx u_{\boldsymbol{\theta}}$ into \cref{eq:data_model} yields the parameter-to-observable map $G(\boldsymbol{\theta}) := F(u_{\boldsymbol{\theta}})$ and the likelihood function
\begin{equation}\label{eq:likelihood}
	f(\boldsymbol{\theta}; \mathbf{y})
	\propto \exp\left( -\frac{1}{2\sigma^2} \norm{ G(\boldsymbol{\theta}) - \mathbf{y} }_2^2 \right).
\end{equation}
For example, if $F \in \R^{M \times N}$ is a linear operator acting on nodal values of $u$ on a grid $\{\mathbf{x}_n\}_{n=1}^{N}$, then $G(\boldsymbol{\theta}) = F \mathbf{u}_{\boldsymbol{\theta}}$ with $\mathbf{u}_{\boldsymbol{\theta}} = [u_{\boldsymbol{\theta}}(\mathbf{x}_1), \dots, u_{\boldsymbol{\theta}}(\mathbf{x}_N)]^\top \in \R^N$.

\begin{remark}[Choice of the activation function]\label{rem:tanh}
	We use $\phi = \tanh$ for the activation function for three reasons.
	(i) 	$\tanh$ is sigmoidal, allowing a single neuron to realize a smoothed jump and making such networks natural parametrizations of piecewise-smooth functions; 
	(ii) $\tanh$ is smooth, rendering $\boldsymbol{\theta} \mapsto G(\boldsymbol{\theta})$ differentiable, as used by the gradient-based solvers in \Cref{sec:MAP}; 
	(iii) the wide-network limit results for heavy-tailed weights that we invoke in \cref{rem:tails} require the activation to be Lipschitz continuous and sub-linearly increasing; see \cite{jung2023alpha,soto2025wide}. 
	The bounded and $1$-Lipschitz $\tanh$ satisfies both conditions trivially.
\end{remark}

\begin{remark}[More general noise models and likelihoods]\label{rem:general_noise}
	The i.i.d.\ assumption on $\mathbf{e}$ is for ease of presentation.
	For zero-mean Gaussian noise with known covariance $\Sigma \in \R^{M \times M}$ the misfit in \cref{eq:likelihood} becomes $|| \Sigma^{-1/2} ( G(\boldsymbol{\theta}) - \mathbf{y} ) ||_2^2$. 
	All subsequent developments apply verbatim.
	Multiplicative noise and combinations of additive and multiplicative components \cite{dunlop2019multiplicative} can be treated similarly.
	In fact, the hierarchical prior construction of \Cref{subsec:GSBL} and the inference schemes of \Cref{sec:MAP,sec:deepEnsembles} interact with the data only through the negative log-likelihood $\boldsymbol{\theta} \mapsto -\log f(\boldsymbol{\theta}; \mathbf{y})$. 
	Any likelihood whose negative logarithm is continuous and bounded from below can therefore be substituted for \cref{eq:likelihood}.
\end{remark}

\subsection{Motivation for sparsity-promoting priors}
\label{subsec:sparse_priors}

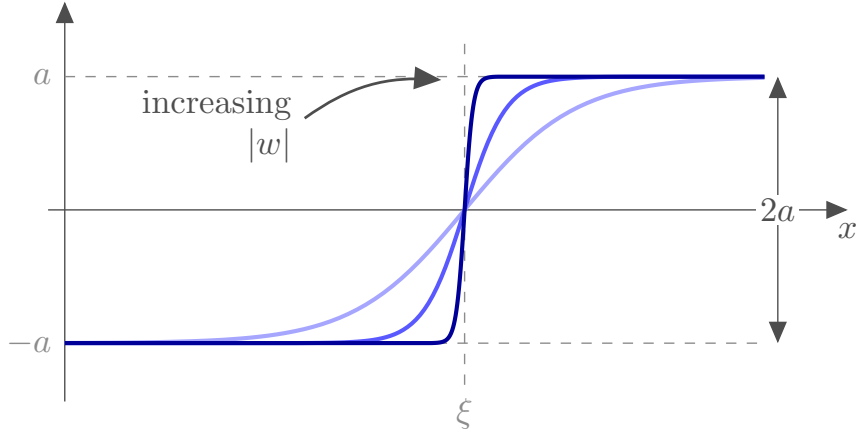
\begin{figure}[tb]
	\centering
	\resizebox{0.7\textwidth}{!}{
\begin{tikzpicture}[x=3.4cm, y=1.7cm, font=\small]
	\pgfmathdeclarefunction{mytanh}{1}{\pgfmathparse{(exp(2*#1)-1)/(exp(2*#1)+1)}}

	\def\aa{0.8}   
	\def\xipos{0.2}   

	\draw[->, black!70] (-1.05, 0) -- (1.35, 0) node[below] {$x$};
	\draw[->, black!70] (-1.0, -1.15) -- (-1.0, 1.25);

	\draw[dashed, black!45] (-1.0, \aa) node[left] {$a$} -- (1.1, \aa);
	\draw[dashed, black!45] (-1.0, -\aa) node[left] {$-a$} -- (1.1, -\aa);
	\draw[dashed, black!45] (\xipos, -1.05) node[below] {$\xi$} -- (\xipos, 1.0);

	\draw[blue!35, very thick]
		plot[domain=-1.0:1.1, samples=80] (\x, {\aa*mytanh(3*(\x-\xipos))});
	\draw[blue!65, very thick]
		(-1.0, {-\aa*mytanh(4)}) -- ({\xipos-0.5}, {-\aa*mytanh(4)})
		plot[domain=\xipos-0.5:\xipos+0.5, samples=80] (\x, {\aa*mytanh(8*(\x-\xipos))})
		-- (1.1, {\aa*mytanh(4)});
	\draw[blue!60!black, very thick]
		(-1.0, {-\aa*mytanh(4)}) -- ({\xipos-0.1}, {-\aa*mytanh(4)})
		plot[domain=\xipos-0.1:\xipos+0.1, samples=80] (\x, {\aa*mytanh(40*(\x-\xipos))})
		-- (1.1, {\aa*mytanh(4)});

	\draw[->, black!70, thick] (-0.28, 0.55) to[bend left=20] (0.13, 0.78);
	\node[black!70, left, align=right] at (-0.28, 0.5) {increasing\\ $\abs{w}$};

	\draw[<->, black!70] (1.14, -\aa) -- (1.14, \aa)
		node[midway, fill=white, inner sep=1pt] {$2a$};
\end{tikzpicture}
	}
	\caption{
		A single neuron $x \mapsto a \tanh\bigl( w (x - \xi) \bigr)$ for fixed amplitude $a$ and location $\xi$, and increasing $\abs{w}$.
		The neuron transitions between $-a$ and $a$ around $\xi$, and the transition sharpens into a jump of height $2a$ as $\abs{w} \to \infty$.
		Forming a sharp transition thus requires a neural parameter $w$ of large magnitude.
	}
	\label{fig:single_neuron}
\end{figure}

The prior on the neural parameters $\boldsymbol{\theta}$ determines whether realizations of $u_{\boldsymbol{\theta}}$ favor smoothness or sharp transitions.
To see what edge-preservation requires, consider a single neuron in one spatial dimension: The mapping $x \mapsto a \tanh\bigl( w (x - \xi) \bigr)$ transitions between the values $-a$ and $a$ around location $\xi$, and the transition sharpens as $\abs{w}$ increases, approaching a jump of height $2a$ as $\abs{w} \to \infty$; see \cref{fig:single_neuron}.
Representing a piecewise-smooth function with a few jumps therefore requires a few neural parameters of large magnitude---producing the sharp transitions---while most parameters should remain small so that $u_{\boldsymbol{\theta}}$ stays smooth away from the jumps.
That is, edge-preservation corresponds to \emph{sparsity in the parameter amplitudes}.

A Gaussian prior penalizes large parameters quadratically and thus suppresses exactly the large-amplitude parameters needed to form sharp transitions.
Heavy-tailed priors, in contrast, reconcile many near-zero parameters with a few large ones. 
This observation underlies the BNN priors of \cite{li2022bayesian,lee2023deep,soto2025wide}.

\subsection{Hierarchical sparsity-promoting priors}
\label{subsec:GSBL}

Instead of directly placing heavy-tailed distributions on the neural parameters, we model them \emph{hierarchically} using conditionally Gaussian scale mixtures \cite{andrews1974scale,west1987scale}.
Specifically, we consider the hierarchical \emph{GSBL prior} model
\begin{equation}\label{eq:GSBL_prior}
	\theta_i \mid \psi_i \sim \mathcal{N}(0, \psi_i),
	\quad
	\psi_i \sim \mathcal{GG}(r, \beta, \vartheta_i),
	\quad i = 1, \dots, K,
\end{equation}
independently, where $\mathcal{GG}(r, \beta, \vartheta)$ denotes the generalized gamma distribution \cite{stacy1962generalization} with density
\begin{equation}\label{eq:GG_density}
	\pi^0(\psi \mid r, \beta, \vartheta)
		= \frac{\abs{r}}{\Gamma(\beta) \vartheta} \left( \frac{\psi}{\vartheta} \right)^{r\beta - 1} \exp\left( - \left( \frac{\psi}{\vartheta} \right)^{r} \right),
		\quad
		\psi > 0,
\end{equation}
parametrized by the parameters $r \in \R \setminus \{0\}$, $\beta > 0$, and $\vartheta > 0$. 
\Cref{fig:pipeline} illustrates the resulting generative process. 
The hyperparameters $\boldsymbol{\psi} = [\psi_1, \dots, \psi_K]^\top$ act as local variances that the inference procedure learns jointly with $\boldsymbol{\theta}$. 
Parameters irrelevant to explaining the data are shrunk aggressively ($\psi_i$ small), while parameters supporting sharp transitions may become large ($\psi_i$ large).
The joint prior, $\pi^0(\boldsymbol{\theta}, \boldsymbol{\psi}) 	= \pi^0(\boldsymbol{\theta} \mid \boldsymbol{\psi}) \, \pi^0(\boldsymbol{\psi})$, decomposes into a conditional Gaussian prior and a generalized gamma hyperprior with densities
\begin{equation}\label{eq:joint_prior}
\resizebox{.9\textwidth}{!}{$\displaystyle 
	\pi^0(\boldsymbol{\theta} \mid \boldsymbol{\psi})
		= \prod_{i=1}^{K} \frac{1}{\sqrt{2\pi \psi_i}} \exp\left( -\frac{\theta_i^2}{2\psi_i} \right), \quad
	\pi^0(\boldsymbol{\psi})
		= \prod_{i=1}^{K} \left( \frac{\abs{r}}{\Gamma(\beta) \vartheta_i} \right) \left( \frac{\psi_i}{\vartheta_i} \right)^{r\beta - 1} \exp\left( - \left( \frac{\psi_i}{\vartheta_i} \right)^{r} \right).
$}
\end{equation}
We call the resulting BNN a \emph{GSBL-BNN}, after the GSBL priors developed for linear inverse problems \cite{calvetti2020sparse,glaubitz2023generalized,lindbloom2025efficient,glaubitz2025efficient}. 
\Cref{tab:special_cases} summarizes important special cases of the GSBL prior \cref{eq:GSBL_prior}, revealing which prior the neural parameters $\boldsymbol{\theta}$ marginally follow when integrating out the hyperparameters $\boldsymbol{\psi}$ in \cref{eq:GSBL_prior}.

\begin{table}[h]
	\centering
	\caption{Special cases of the GSBL prior \cref{eq:GSBL_prior} and the induced marginal distributions of $\theta_i$.}
	\label{tab:special_cases}
	\begin{tabular}{llll}
		\toprule
		$(r, \beta)$ & Hyperprior on $\psi_i$ & Marginal of $\theta_i$ & Tail behavior \\
		\midrule
		$(1, \beta)$ & gamma & variance-gamma & $\log$-linear (Laplace-type) \\
		$(1, 1)$ & exponential & Laplace$\bigl(0, \sqrt{\vartheta_i/2}\bigr)$ & $\log$-linear \\
		$(-1, \beta)$ & inverse gamma & Student-$t$, $2\beta$ d.o.f. & polynomial, index $2\beta$ \\
		$(-1, 1/2)$ & inverse gamma & Cauchy$\bigl(0, \sqrt{2\vartheta_i}\bigr)$ & polynomial, index $1$ \\
		\bottomrule
	\end{tabular}
\end{table}

Before proceeding, let us motivate modeling the neural parameters hierarchically via \cref{eq:GSBL_prior} rather than directly equipping them with the corresponding marginal heavy-tailed priors, as in \cite{li2022bayesian}.
The reasons are computational.
First, conditioned on the hyperparameters $\boldsymbol{\psi}$, the prior on the neural parameters is Gaussian, so that updating $\boldsymbol{\theta}$ reduces to a smooth, quadratically regularized training problem. 
At the same time, conditioned on $\boldsymbol{\theta}$, updating the hyperparameters $\boldsymbol{\psi}$ decouples into $K$ scalar problems that can be solved exactly.
Directly optimizing under marginal heavy-tailed priors instead faces log-densities whose gradients vanish in the tails, yielding flat optimization landscapes. 
Moreover, closed-form expressions for the marginals are available only in special cases. 
Second, the hierarchical formulation is exact rather than approximate. 
In contrast, $\alpha$-stable priors rely on numerical density approximations \cite{suuronen2023bayesian}.

\begin{remark}[Notation]\label{rem:notation}
	We deviate from the GSBL literature convention, where $\boldsymbol{x}$ denotes the unknown and $\boldsymbol{\theta}$ the hyperparameters.
	To remain consistent with the established BNN convention, we use $\boldsymbol{\theta}$ for the neural parameters and $\boldsymbol{\psi}$ for the hyperparameters.
\end{remark}

\begin{remark}[Tail behavior of the marginal priors]\label{rem:tails}
	The tail column of \Cref{tab:special_cases} reflects a general pattern.
	Write $\theta_i = \sqrt{\psi_i}\, Z_i$ with $Z_i \sim \mathcal{N}(0, 1)$ independent of $\psi_i$.
	Since $(\psi_i/\vartheta_i)^{r}$ is gamma distributed with shape $\beta$ and unit scale, the tail probability $\mathbb{P}(\psi_i > \tau)$ decays like $\tau^{-\abs{r}\beta}$ as $\tau \to \infty$ for $r < 0$.
	Breiman's lemma \cite{breiman1965some} transfers this power law to the marginal: $\mathbb{P}(\abs{\theta_i} > t)$ decays like $t^{-\alpha}$ as $t \to \infty$, with \emph{tail index} $\alpha = 2\abs{r}\beta$.
	Hence, $\mathbb{E}\abs{\theta_i}^p < \infty$ if and only if $p < \alpha$, and $\theta_i$ has infinite variance if and only if $\alpha \leq 2$.
	For $\alpha < 2$, the marginals lie in the domain of attraction of a symmetric $\alpha$-stable law \cite{gnedenko1968limit,samorodnitsky1994stable}. 
	The GSBL-BNN prior thus inherits the edge-preserving behavior of the heavy-tailed BNN priors of \cite{li2022bayesian}, whose suitably scaled wide-network limits are $\alpha$-stable processes with discontinuity-like realizations \cite{soto2025wide,lee2023deep}.
	For $r > 0$, in contrast, $\mathbb{P}(\psi_i > \tau)$ decays faster than any power of $\tau$. 
	Hence, all moments of $\theta_i$ are finite and wide-network limits are Gaussian-process-like \cite{neal1996bayesian} and not edge-preserving.
\end{remark}

\subsection{Specifying the hierarchical prior}
\label{subsec:prior_spec}

We comment on a few remaining modeling choices that we found to be important for the performance of the proposed GSBL-BNN priors.

\subsubsection{The role of the weights and biases}
\label{subsubsec:prior_spec_roles}

The weights and biases play distinct roles, as identified in \Cref{subsec:sparse_priors}. 
The \emph{weights control the amplitude and sharpness of the transitions} realized by the individual neurons---which is where edge-preservation calls for the sparsity of the GSBL prior---whereas the \emph{biases position those transitions}.
In the single-neuron example of \Cref{subsec:sparse_priors}, a first-layer neuron transitions at $\xi = -b/w$, where its argument $w x + b$ vanishes.
Sparsity thus offers no benefit for the biases. Shrinking a bias toward zero would pin its transition to a fixed reference point rather than let it move freely across the domain, which was also noted in \cite{li2022bayesian}.
The biases therefore require a different treatment, which we determine next by asking where the induced transitions should lie.

\subsubsection{Choosing the bias distribution}
\label{subsubsec:prior_spec_bias}

\begin{figure}[tb]
	\centering
	\includegraphics[width=\textwidth]{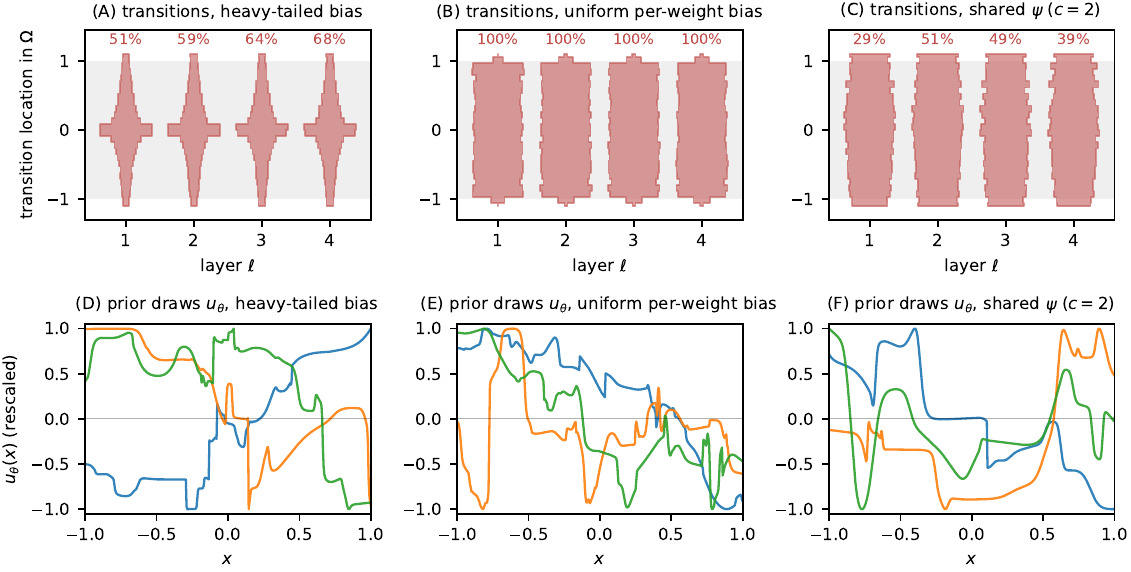}
	\caption{
		Three bias schemes for the GSBL-BNN prior with $r = -1$, $\beta = 1/2$, $\vartheta^{(1)} = 100$, and $\vartheta^{(\ell)} = 64^{-2}$ with tail index $\alpha = 1$ in all deeper layers, including the output layer; see \cref{eq:fan_in}.
	\emph{(A)--(C):} Per-layer distribution of the transition locations, where a neuron's pre-activation changes sign, and percentage of neurons transitioning within $\Omega = [-1, 1]$, over $400$ random networks with four hidden layers of width $64$.
	\emph{(D)--(F):} Prior draws of $u_{\boldsymbol{\theta}}$, each shifted and scaled to span $[-1, 1]$.
	}
	\label{fig:bias_calibration}
\end{figure}

For an edge-preserving prior, the transitions should be free to lie anywhere in $\Omega$ rather than cluster.
The obstacle is that a neuron's transition location $\xi = -b/w$ is a \emph{ratio}; what fixes the edge is the size of the bias \emph{relative} to the weight.
A heavy-tailed weight prior makes $\abs{w}$ vary wildly from neuron to neuron, so a bias whose scale cannot track each neuron's own $\abs{w}$ yields a small ratio $\abs{b/w}$ precisely for the neurons with large weights---driving their transitions toward $\xi = 0$ and piling the jumps of $u_{\boldsymbol{\theta}}$ near the center of $\Omega$.
Placing edges freely thus requires \emph{coupling} the bias to its weight.
We discuss three options for treating the biases, which are illustrated in \Cref{fig:bias_calibration}.

\paragraph{Independent biases (existing edge-preserving BNNs)}

The heavy-tailed BNNs of \cite{li2022bayesian} draw the biases independently of the weights---typically from a heavy-tailed distribution of their own.
As a result, the edges concentrate near the center; see \Cref{fig:bias_calibration}(A, D).

\paragraph{Uniform per-weight biases}

Drawing each bias uniformly over the range of pre-activations its neuron attains on $\Omega$ does spread the edges. 
Each neuron then transitions at a uniformly random point of the domain, as illustrated in \Cref{fig:bias_calibration}(B). 
Still, this approach has two drawbacks: 
It activates every neuron, giving `busy' realizations rather than ones with a few distinct jumps; see \Cref{fig:bias_calibration}(E).
Second, the induced bias prior is uniform on a weight-dependent interval, a hard constraint that is not conditionally Gaussian and thus breaks the structure that makes the MAP estimation of \Cref{sec:MAP} efficient.

\paragraph{Shared-$\psi$ biases (ours)}

\begin{figure}[tb]
	\centering
\begin{tikzpicture}[
	neuron/.style={circle, draw, minimum size=0.95cm, inner sep=1pt, font=\small},
	hiddenA/.style={neuron, draw=blue!55!black, fill=blue!10},
	hiddenB/.style={neuron, draw=orange!75!black, fill=yellow!22},
	biasnode/.style={neuron, draw=green!55!black, fill=green!12},
	psinode/.style={neuron, draw=red!65!black, fill=red!8, dashed},
	conn/.style={->, semithick, black!65},
	hyper/.style={->, dashed, semithick, red!65!black, shorten >=1pt},
	lab/.style={font=\small}
]
	\node[hiddenA] (h1) at (0, 1.8) {$h_1^{(\ell-1)}$};
	\node[hiddenA] (h2) at (0, 0.6) {$h_2^{(\ell-1)}$};
	\node (hdots) at (0, -0.15) {$\vdots$};
	\node[hiddenA] (hn) at (0, -1.1) {$h_{n_{\ell-1}}^{(\ell-1)}$};

	\node[hiddenB] (hi) at (5.2, 0) {$h_i^{(\ell)}$};
	\draw[conn] (hi) -- ++(1.0, 0);
	\node[biasnode] (one) at (5.2, 1.8) {$1$};
	\draw[conn, green!55!black] (one) -- node[lab, right, pos=0.5, green!55!black] {$b_i^{(\ell)} \sim \mathcal{N}\bigl(0, c^2 \psi_i^{(\ell)}\bigr)$} (hi);

	\draw[conn] (h1) -- node[lab, above, sloped, pos=0.25, blue!55!black] {$W_{i,1}^{(\ell)}$} (hi);
	\draw[conn] (h2) -- node[lab, above, sloped, pos=0.25, blue!55!black] {$W_{i,2}^{(\ell)}$} (hi);
	\draw[conn] (hn) -- node[lab, above, sloped, pos=0.25, blue!55!black] {$W_{i,n_{\ell-1}}^{(\ell)}$} (hi);

	\node[psinode] (psi) at (2.6, 2.6) {$\psi_i^{(\ell)}$};
	\node[lab, red!65!black, anchor=east] at (2.05, 2.7) {$\psi_i^{(\ell)} \sim \mathcal{GG}\bigl(r, \beta, \vartheta^{(\ell)}\bigr)$};
	\node[lab, blue!55!black, anchor=east] at (-0.65, 0.6) {$\mathbf{w}_i^{(\ell)} \sim \mathcal{N}\bigl(\mathbf{0}, \psi_i^{(\ell)} I\bigr)$};
	\draw[hyper] (psi) -- (2.6, 1.05);
	\draw[hyper] (psi) to[bend right=10] (5.07, 1.15);
\end{tikzpicture}
	\caption{
		Shared-$\psi$ prior for the parameters of a single hidden neuron $h_i^{(\ell)}$ of the network in \cref{fig:NN_schematic}. 
		Its incoming weights and bias form one parameter block $\boldsymbol{\theta}_i^{(\ell)} = (\boldsymbol{w}_i^{(\ell)}, b_i^{(\ell)})$, governed by a single shared hyperparameter $\psi_i^{(\ell)}$.
	}
	\label{fig:shared_psi}
\end{figure}
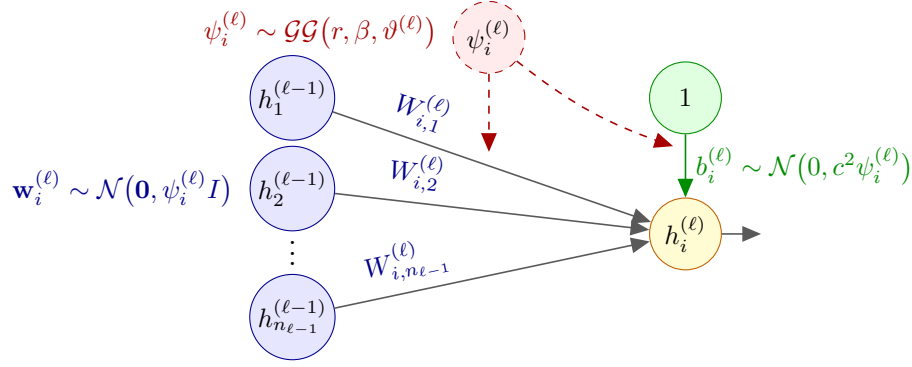

We couple each hidden neuron's bias to its weights through a shared hyperparameter. 
Specifically, a neuron's incoming weights are $\mathcal{N}(0, \psi)$ and its bias is $\mathcal{N}(0, c^2 \psi)$, sharing a common variance $\psi \sim \mathcal{GG}(r, \beta, \vartheta)$ with a tunable constant $c > 0$.
\Cref{fig:shared_psi} illustrates this coupling.
For neuron $i$ in layer $\ell$, the $n_{\ell-1}$ incoming weights $\boldsymbol{w}_i^{(\ell)}$---the $i$th row of $W^{(\ell)}$---and the bias $b_i^{(\ell)}$ form one parameter block $\boldsymbol{\theta}_i^{(\ell)} = (\boldsymbol{w}_i^{(\ell)}, b_i^{(\ell)}) \in \R^{n_{\ell-1}+1}$ governed by a single shared hyperparameter $\psi_i^{(\ell)}$.
Writing $b = c\sqrt{\psi}\, Z_b$ and $w = \sqrt{\psi}\, Z_w$ with independent standard normals $Z_b, Z_w \sim \mathcal{N}(0, 1)$, the transition location is then $\xi = -b/w = -c\, Z_b/Z_w$, a scaled ratio of standard normals.
Notably, $\xi$ thus follows a scaled Cauchy distribution $\mathrm{Cauchy}(0, c)$, \emph{independent} of $\psi$ and of the hyperprior. 
Consequently, the edges spread roughly evenly across $\Omega$ for every $(r, \beta)$, which we illustrate in \Cref{fig:bias_calibration}(C, F).

To our knowledge, this coupling is new among edge-preserving BNNs, and it removes both previous shortcomings: 
(i) Unlike the independent biases of \cite{li2022bayesian}, our shared-$\psi$ approach frees the edge locations and does not cluster them close to the domain's center. 
(ii) Unlike the uniform bias, our approach keeps the biases conditionally Gaussian, preserving the structure we leverage for efficient MAP estimation in \Cref{sec:MAP}.

The scale $c$ is a free parameter, tuned alongside $(r, \beta, \vartheta)$. 
It trades evenness of edge placement against neuron activity.
The density of $\xi$ over $\Omega$ has center-to-edge ratio $1 + 1/c^2$, flattening toward uniform as $c$ grows, while the fraction of neurons transitioning within $\Omega$ is $(2/\pi)\arctan(1/c)$.
In \Cref{fig:bias_calibration}(F) we take $c = 2$ as a representative choice.

\subsubsection{Calibration of the scales}
\label{subsubsec:prior_spec_vartheta}

The scale $\vartheta_i$ sets the magnitude of $\theta_i$ in the following sense.
Because the generalized gamma density \cref{eq:GG_density} depends on $\psi$ only through $\psi/\vartheta$, the hyperparameter factors as $\psi_i = \vartheta_i\, V_i$ with random variable $V_i \sim \mathcal{GG}(r, \beta, 1)$.
Hence, $\theta_i = \sqrt{\psi_i}\, Z_i = \sqrt{\vartheta_i}\, (\sqrt{V_i \, Z_i} )$ with $Z_i \sim \mathcal{N}(0, 1)$ as in \cref{rem:tails}.
Every scale-type summary of $|\theta_i|$, such as its median and its standard deviation, is therefore proportional to $\sqrt{\vartheta_i}$. 

We choose the scales per layer so that prior realizations of $u_{\boldsymbol{\theta}}$ remain order one as the widths grow, decoupling the shape parameters $(r, \beta)$ of the GSBL prior from the architecture.
The relevant quantity is a hidden pre-activation:
\begin{equation}\label{eq:preactivation}
	z_i^{(\ell)} = \sum_{j=1}^{n_{\ell-1}} W_{ij}^{(\ell)} h_j^{(\ell-1)} + b_i^{(\ell)},
\end{equation}	
which is a sum of $n_{\ell-1}$ independent weight contributions of scale $\sqrt{\vartheta^{(\ell)}}$ weighted by the bounded activations $h_j^{(\ell-1)} \in [-1, 1]$.
The weights have tail index $\alpha = 2\abs{r}\beta$ (see \cref{rem:tails}), and a sum of $n$ independent terms of unit scale and tail index $\alpha$ has scale $\sim n^{1/\min(\alpha, 2)}$. 
Notably, this scale becomes the central-limit rate $n^{1/2}$ when the variance is finite ($\alpha \ge 2$) and the $\alpha$-stable rate $n^{1/\alpha}$ when it is infinite ($\alpha < 2$); see \cite{gnedenko1968limit,samorodnitsky1994stable}.
Keeping $z_i^{(\ell)}$ of order one as $n_{\ell-1} \to \infty$ thus requires $\sqrt{\vartheta^{(\ell)}} \sim n_{\ell-1}^{-1/\min(\alpha, 2)}$. 
We thus choose 
\begin{equation}\label{eq:fan_in}
	\vartheta^{(\ell)} = 
	\begin{cases} 
		\kappa n_{\ell-1}^{-1} & \text{if } \alpha \ge 2, \\ 
		\kappa n_{\ell-1}^{-2/\alpha} & \text{if } \alpha < 2,
	\end{cases}
\end{equation}
for all layers except the first, whose scale we treat separately below.
Here, $\kappa > 0$ is a width-independent constant that controls how strongly the deeper hidden neurons saturate and, through the linear output layer, the amplitude of $u_{\boldsymbol{\theta}}$.

\begin{figure}[tb]
	\centering
	\includegraphics[width=\textwidth]{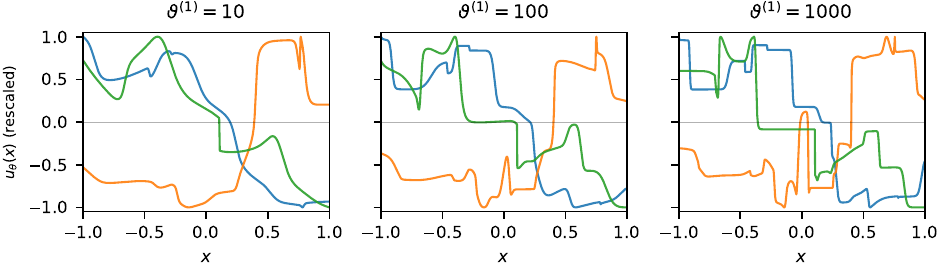}
	\caption{
		Prior draws of $u_{\boldsymbol{\theta}}$ for increasing first-layer scale $\vartheta^{(1)}$, using the shared-$\psi$ prior and networks of \Cref{fig:bias_calibration} ($r = -1$, $\beta = 1/2$, $c = 2$).
	Each panel shows three draws, each shifted and scaled to span $[-1, 1]$.
	As $\vartheta^{(1)}$ grows, the jumps sharpen while largely staying in place, since the first-layer transition locations $\xi = -b/w$ do not depend on $\vartheta^{(1)}$.
	}
	\label{fig:theta1_sweep}
\end{figure}

The first-layer scale $\vartheta^{(1)}$ additionally controls the \emph{sharpness} of the induced jumps.
A first-layer neuron $a\tanh(w(x - \xi))$ swings between its saturated values $\pm a$ while its argument traverses an $\mathcal{O}(1)$ window. 
Notably, $\tanh(\pm 2) \approx \pm 0.96$ is already near saturation, so its transition occupies an input width $\abs{x - \xi} \lesssim 2/\abs{w}$.
Since $\abs{w} \sim \sqrt{\vartheta^{(1)}}$ by the scale argument above, the input width is $\abs{x - \xi} \sim (\vartheta^{(1)})^{-1/2}$. 
Thus, enlarging $\vartheta^{(1)}$ steepens the transitions, and, because $\tanh$ is bounded, does so without inflating the amplitude of $u_{\boldsymbol{\theta}}$. 
At the same time, a common rescaling of all $\vartheta^{(\ell)}$ would instead inflate the amplitude through the linear output layer.
We therefore treat $\vartheta^{(1)}$ as another free parameter governing feature sharpness, illustrated in \Cref{fig:theta1_sweep}; the draws in \Cref{fig:bias_calibration} use $\vartheta^{(1)} = 100$.
Both figures use $\kappa = 1$.

\subsection{The joint posterior}
\label{subsec:posterior}

Combining the likelihood \cref{eq:likelihood} with the hierarchical prior \cref{eq:joint_prior}, Bayes' theorem yields the joint posterior density
\begin{equation}\label{eq:joint_posterior}
	\pi^{\mathbf{y}}(\boldsymbol{\theta}, \boldsymbol{\psi})
	\propto f(\boldsymbol{\theta}; \mathbf{y})\, \pi^0(\boldsymbol{\theta} \mid \boldsymbol{\psi})\, \pi^0(\boldsymbol{\psi})
\end{equation}
over the neural parameters \emph{and} the hyperparameters.
Writing $D_{\boldsymbol{\psi}} := \diag(\boldsymbol{\psi})$ and substituting the likelihood \cref{eq:likelihood} and the GSBL prior \cref{eq:joint_prior}, we get 
\begin{equation}\label{eq:joint_posterior_explicit}
\resizebox{.9\textwidth}{!}{$\displaystyle 
	\pi^{\mathbf{y}}(\boldsymbol{\theta}, \boldsymbol{\psi})
		\propto \exp\left( 
			- \frac{1}{2\sigma^2} \| G(\boldsymbol{\theta}) - \mathbf{y} \|_2^2 
			- \frac{1}{2} \| D_{\boldsymbol{\psi}}^{-1/2} \boldsymbol{\theta} \|_2^2 
			+ \sum_{i=1}^K \left[ 
				\left( r \beta - \frac{3}{2} \right) \log \psi_i 
				- \left( \frac{\psi_i}{\vartheta_i} \right)^{r} 
			\right]
		\right)
$}
\end{equation}
for the joint posterior \cref{eq:joint_posterior}.
\section{Efficient MAP estimation via block-coordinate descent}
\label{sec:MAP}

We now discuss efficient and robust MAP estimation for the GSBL-BNN posterior \cref{eq:joint_posterior} using a BCD approach that alternates between updating the neural parameters and hyperparameters.

\subsection{The MAP objective}
\label{subsec:objective}

Taking the negative logarithm of \cref{eq:joint_posterior_explicit} and discarding additive constants, the MAP estimates of \cref{eq:joint_posterior_explicit} are the minimizers of the objective
\begin{equation}\label{eq:objective}
	J(\boldsymbol{\theta}, \boldsymbol{\psi})
	:= \frac{1}{2\sigma^2} \| G(\boldsymbol{\theta}) - \mathbf{y} \|_2^2 
		+ \frac{1}{2} \| D_{\boldsymbol{\psi}}^{-1/2} \boldsymbol{\theta} \|_2^2 
		+ \sum_{i=1}^K \left[ 
				\left( \frac{\psi_i}{\vartheta_i} \right)^{r} 
				- \left( r \beta - \frac{3}{2} \right) \log \psi_i
			\right]
\end{equation}
over $(\boldsymbol{\theta}, \boldsymbol{\psi}) \in \R^K \times \R^K_{>0}$.
That is,
\begin{equation}\label{eq:MAP_problem}
	(\boldsymbol{\theta}_{\MAP}, \boldsymbol{\psi}_{\MAP})
	\in \argmin_{\boldsymbol{\theta} \in \R^K,\, \boldsymbol{\psi} \in \R^K_{>0}} J(\boldsymbol{\theta}, \boldsymbol{\psi}).
\end{equation}

Throughout, we assume that $G$ is continuously differentiable, which holds for the network \cref{eq:NN} whenever $F$ is, and that $r < 0$ or $r\beta > 3/2$, since otherwise $J$ may have no minimizer.
For the heavy-tailed regime $r < 0$ advocated in \cref{rem:tails}, $J$ is moreover coercive, i.e., $J(\boldsymbol{\theta}, \boldsymbol{\psi}) \to \infty$ whenever $\norm{(\boldsymbol{\theta}, \boldsymbol{\psi})} \to \infty$ or $\min_i \psi_i \to 0$, so that minimizers exist.




\subsection{The BCD algorithm}
\label{subsec:BCD}

Inspired by efficient MAP estimation for linear inverse problems with hierarchical priors \cite{calvetti2020sparse,glaubitz2023generalized,dong2023inducing,xiao2023sequential,glaubitz2024leveraging,lindbloom2025efficient,lindbloom2025priorconditioned,glaubitz2025bayesian}, we propose the following BCD approach to minimize $J$ by alternating between the two blocks $\boldsymbol{\theta}$ and $\boldsymbol{\psi}$:
\begin{equation}\label{eq:BCD}
	\boldsymbol{\theta}^{(k+1)} 
		\approx \argmin_{\boldsymbol{\theta} \in \R^K} J(\boldsymbol{\theta}, \boldsymbol{\psi}^{(k)}),
	\quad
	\boldsymbol{\psi}^{(k+1)} 
		= \argmin_{\boldsymbol{\psi} \in \R^K_{>0}} J(\boldsymbol{\theta}^{(k+1)}, \boldsymbol{\psi}),
	\quad k = 0, 1, 2, \dots
\end{equation}
The motivation for the above BCD approach is that the two subproblems in \cref{eq:BCD} are lower-dimensional and easier to solve than the joint optimization problem \cref{eq:MAP_problem}.

\subsubsection*{The $\boldsymbol{\theta}$-update} 

The $\boldsymbol{\theta}$-update in \cref{eq:BCD} is standard network training with iteratively reweighted $\ell^2$-regularization.
For fixed $\boldsymbol{\psi} = \boldsymbol{\psi}^{(k)}$, all terms of \cref{eq:objective} that depend on $\boldsymbol{\theta}$ reduce to
\begin{equation}\label{eq:w_step}
	J(\boldsymbol{\theta}, \boldsymbol{\psi}^{(k)})
		= \frac{1}{2\sigma^2} \norm{ G(\boldsymbol{\theta}) - \mathbf{y} }_2^2
		+ \frac{1}{2} \norm{ D_{\boldsymbol{\psi}^{(k)}}^{-1/2} \boldsymbol{\theta} }_2^2
\end{equation}
with $D_{\boldsymbol{\psi}^{(k)}} = \diag( \boldsymbol{\psi}^{(k)} )$.
The first term on the right-hand side of \cref{eq:w_step} is the standard data-misfit objective of network training, and the second is a weighted $\ell^2$-penalty that shrinks each parameter in proportion to $1/\psi_i^{(k)}$.
The penalty weights $1/\psi_i^{(k)}$ are held fixed within the $\boldsymbol{\theta}$-update but are recomputed by the $\boldsymbol{\psi}$-update at every outer iteration. 

Any standard first-order optimizer can be used to minimize \cref{eq:w_step} w.\,r.\,t.\ $\boldsymbol{\theta}$.
In \Cref{sec:numerics}, we use Muon without weight decay~\cite{jordan2024muon} with warm starts, initializing the inner iteration at the previous iterate $\boldsymbol{\theta}^{(k)}$.
Using Adam instead of Muon, we obtained similar results, which we omit for brevity.
We perform refinement with the limited memory Broyden-Fletcher-Goldfarb-Shanno (L-BFGS) optimizer~\cite{liu_limited_1989}.

\subsubsection*{The $\boldsymbol{\psi}$-update}

The $\boldsymbol{\psi}$-update in \cref{eq:BCD} can be computed exactly.
We first provide the $\boldsymbol{\psi}$-update for the case in which the biases independently follow the same GSBL prior as the weights (see the first paragraph in \Cref{subsubsec:prior_spec_bias}). 
We explain how the $\boldsymbol{\psi}$-update changes for the shared-$\psi$ biases case (see the third paragraph in \Cref{subsubsec:prior_spec_bias}) in \cref{rem:shared}.

For fixed $\boldsymbol{\theta} = \boldsymbol{\theta}^{(k+1)}$, the objective \cref{eq:objective} decouples into $K$ independent scalar problems:
\begin{equation}\label{eq:theta_step}
	\psi_i^{(k+1)}
		= \argmin_{\psi > 0}\left\{ g_{\theta_i^{(k+1)}}(\psi) \right\}
		\quad \text{with} \quad
	g_{\theta}(\psi)
	:= \frac{\theta^2}{2\psi} + \left( \frac{\psi}{\vartheta} \right)^{r} - \eta \log \psi.
\end{equation}
We write $\eta = r \beta - 3/2$ and $\vartheta$ for $\vartheta_i$ to simplify notation.
Notably, for $r = \pm 1$, \cref{eq:theta_step} admits a closed-form solution: 
Setting the derivative of $g_{\theta_i^{(k+1)}}(\psi)$ to zero and solving it for $\psi$ yields 
\begin{equation}\label{eq:theta_update_closed}
	\psi^{\ast}(\theta)
		= \frac{\vartheta}{2} \left( \eta + \sqrt{ \eta^2 + \frac{2 \theta^2}{\vartheta} } \right) \ \text{ if } r = 1, \quad 
	\psi^{\ast}(\theta)
		= \frac{\vartheta + \theta^2/2}{\beta + 3/2} \ \text{ if } r = -1.
\end{equation}
For general $r$, the root of $g_{\theta_i^{(k+1)}}'(\psi)$ can be computed by a safeguarded Newton iteration.
Alternatively, we can obtain all updates by numerically solving a single scalar ordinary differential equation (ODE), as in \cite{calvetti2020sparse,lindbloom2025efficient}.

\begin{remark}[Modifications for shared-$\psi$ biases]\label{rem:shared}
	Under our proposed bias--weight coupling of \Cref{subsec:prior_spec}, the weights are $\mathcal{N}(\mathbf{0}, \psi_i^{(\ell)} I)$ while the bias is $\mathcal{N}(0, c^2 \psi_i^{(\ell)})$. 
	Writing $\| \boldsymbol{\theta}_i^{(\ell)} \|_{c}^2 := \| \boldsymbol{w}_i^{(\ell)} \|_2^2 + (b_i^{(\ell)})^2/c^2$ and $\eta^{(\ell)} = r\beta - 1 - (n_{\ell-1} + 1)/2$, the objective \cref{eq:objective} thus becomes
	\begin{equation}\label{eq:objective_tied}
		J(\boldsymbol{\theta}, \boldsymbol{\psi})
		= \frac{1}{2\sigma^2} \| G(\boldsymbol{\theta}) - \mathbf{y} \|_2^2
		+ \sum_{\ell=1}^{L} \sum_{i=1}^{n_\ell} \left[
			\frac{\| \boldsymbol{\theta}_i^{(\ell)} \|_{c}^2}{2 \psi_i^{(\ell)}}
			+ \left( \frac{\psi_i^{(\ell)}}{\vartheta^{(\ell)}} \right)^{r}
			- \eta^{(\ell)} \log \psi_i^{(\ell)}
		\right].
	\end{equation}
	For fixed $\boldsymbol{\theta}$, the sum in \cref{eq:objective_tied} decouples across blocks, and each block subproblem has exactly the scalar form \cref{eq:theta_step}.
	Hence the solution of \cref{eq:theta_step} carries over to each block verbatim, with $\theta^2$ replaced by $\| \boldsymbol{\theta}_i^{(\ell)} \|_{c}^2$ and $\eta$ by $\eta^{(\ell)}$.
	Specifically, for $r = \pm 1$, we get the closed-form solutions  
\begin{equation}
\resizebox{.9\textwidth}{!}{$\displaystyle 
	\psi_i^{(\ell)\ast}
		= \frac{\vartheta^{(\ell)}}{2} \left( \eta^{(\ell)} + \sqrt{ (\eta^{(\ell)})^2 + \frac{2 \| \boldsymbol{\theta}_i^{(\ell)} \|_{c}^2}{\vartheta^{(\ell)}} } \right) 
		\text{if } r = 1, \quad 
	\psi_i^{(\ell)\ast}
		= \frac{\vartheta^{(\ell)} + \| \boldsymbol{\theta}_i^{(\ell)} \|_{c}^2 / 2 }{ \beta + 1 + (n_{\ell-1} + 1)/2} 
		\text{ if } r = -1.
$}
\end{equation}
	For general $r$, we again use safeguarded Newton iterations or the ODE approach  \cite{calvetti2020sparse,lindbloom2025efficient}.
\end{remark}

We summarize the resulting BCD method in \cref{alg:BCD}. 
We discuss how we initialize training for GSBL-BNNs in \cref{rem:ensemble_init}, where we consider a whole ensemble of networks.

\begin{algorithm}[tb]
	\caption{BCD approach for MAP estimation of GSBL-BNNs}
	\label{alg:BCD}
	\begin{algorithmic}[1]
		\Require initialization $(\boldsymbol{\theta}^{(0)}, \boldsymbol{\psi}^{(0)})$
		\For{$k = 0, 1, 2, \dots$ until converged}
			\State $\boldsymbol{\theta}^{(k+1)} \gets$ approximate minimizer of \cref{eq:w_step}, warm-started at $\boldsymbol{\theta}^{(k)}$
			\State $\psi_i^{(k+1)} \gets \psi^{\ast}(\theta_i^{(k+1)})$ for $i = 1, \dots, K$ \Comment{closed form, Newton, or ODE}
		\EndFor
		\State \Return $(\boldsymbol{\theta}^{(k+1)}, \boldsymbol{\psi}^{(k+1)})$
	\end{algorithmic}
\end{algorithm} 
\section{Deep ensembles}
\label{sec:deepEnsembles}

A single run of \cref{alg:BCD} produces one MAP estimate and therefore no UQ.
This section turns the algorithm into an approximate Bayesian inference method by running it as a \emph{deep ensemble} \cite{lakshminarayanan2017simple}. 
We compute $S$ MAP estimates from independent random initializations, and their spread serves as an approximate posterior uncertainty.

The GSBL-BNN posterior \cref{eq:joint_posterior} is challenging to sample: the nonconvexity of the data misfit typically produces genuinely distinct modes corresponding to different plausible reconstructions.
Deep ensembles turn this limitation into a strategy. 
Because independent initializations tend to fall into distinct basins of attraction, the members capture separate posterior modes that single-mode approximations such as the Laplace or mean-field variational methods miss \cite{fort2019deep,izmailov2021bayesian}. 
Concretely, each member is an independent, embarrassingly parallel run of \cref{alg:BCD}, started from a random initialization; see \cref{rem:ensemble_init}.
Denoting the resulting MAP estimates by $\{ (\boldsymbol{\theta}_{(s)}, \boldsymbol{\psi}_{(s)}) \}_{s=1}^{S}$ and the induced reconstructions by $u_{(s)} := u_{\boldsymbol{\theta}_{(s)}}$, we treat $\{ u_{(s)} \}_{s=1}^{S}$ as approximate samples from the predictive posterior and estimate posterior expectations of functionals of $u$ by the corresponding ensemble averages.

We emphasize that deep ensembles are a heuristic. The ensemble members are stationary points of the objective \cref{eq:objective} rather than draws from \cref{eq:joint_posterior}, and their spread conflates posterior structure with optimization randomness.
We nevertheless adopt them for their simplicity.
The accuracy and calibration of the resulting estimates are assessed in \Cref{sec:numerics}.

\begin{remark}[Initialization strategies]\label{rem:ensemble_init}
	The randomization in MAP estimation (\cref{alg:BCD}) enters through the initialization, for which several strategies are available.
	(i) \emph{Prior draws}, used in \cite{glaubitz2025efficient}, sample $(\boldsymbol{\psi}^{0}, \boldsymbol{\theta}^{0})$ from the hierarchical prior \cref{eq:GSBL_prior}.
	This is arguably the natural choice, producing model-consistent starting points as diverse as the prior itself.
	(ii) \emph{Standard network initializations}, such as LeCun, Glorot, or He scaling~\cite{lecun_efficient_1998,glorot2010understanding,he2015delving}, draw $\boldsymbol{\theta}^{0}$ from distributions whose variance is scaled by the layer fan-in (and, for Glorot, also the fan-out). 
	We then obtain $\boldsymbol{\psi}^{0}$ from a single $\boldsymbol{\psi}$-update. 
	These initializations are tuned for trainability rather than prior consistency.
\end{remark}

\begin{remark}[Alternative ensemble strategies]\label{rem:ensemble_alternatives}
	Random-restart ensembles, considered above, sit at the simple end of a spectrum of increasingly elaborate posterior approximations. 
	Since all members minimize the same objective, an effective optimizer can drive them together and cause underdispersion.
	Randomized-MAP (or randomize-then-optimize) ensembles \cite{bardsley2012mcmc,bardsley2014randomize,pearce2020uncertainty} address this by perturbing each member's data and anchoring it at a random prior draw, so that the spread of the minimizers mimics the posterior---exactly so in the linear-Gaussian case.
	Local Gaussian approximations, such as the Laplace approximation \cite{daxberger2021laplace} or SWAG \cite{maddox2019simple}, instead fit a Gaussian around a single mode.
	Multi-basin mixtures such as MultiSWAG \cite{wilson2020bayesian} combine the two.
	Full approximate inference via variational methods \cite{blundell2015weight} or stochastic-gradient MCMC \cite{welling2011bayesian} targets the entire posterior at substantially higher cost.
	A systematic comparison of these schemes is beyond the scope of this paper.
\end{remark}
\section{Computational examples}
\label{sec:numerics}

We test the proposed GSBL-BNNs for three different inverse problems: denoising a piecewise constant signal, deblurring a piecewise smooth signal, and inverting a two-dimensional hyperbolic conservation law.
See~\cite{jan_habscheid_2026_23184383} for our openly available reproducibility repository.

\subsection{Implementation and diagnostics}
\label{subsec:num_setup}

We implement all computational examples in JAX~\cite{jax2018github} for automatic differentiation and use Equinox~\cite{kidger2021equinox} to build the neural networks.
For the neural network training, we use the optimization library Optax~\cite{deepmind2020jax_optax_citation}.
We use the Muon optimizer\footnote{We obtained similar results using Adam, which we omit for brevity.} without weight decay~\cite{jordan2024muon} for initial training, followed by L-BFGS~\cite{liu_limited_1989} for refinement.
We use a fixed learning rate of $10^{-3}$, as the loss landscape for the network parameters $\theta$ is transformed every 10 inner iterations after the $\boldsymbol{\psi}$-update is applied.
As discussed in~\cref{rem:tanh}, we use the tanh activation function.

We use the following diagnostics to assess reconstruction accuracy and the quality of uncertainty estimates.
Throughout, we evaluate all reconstructions on a uniform reference grid on $\Omega$ and identify functions with their nodal vectors, denoted by upright bold letters.
In particular, $\mathbf{u}^{\dagger}$ denotes the nodal vector of the ground truth $u^{\dagger}$ and $\mathbf{u}_{(s)}$ that of the $s$th ensemble member $u_{(s)}$.
We summarize each ensemble by the pointwise mean $\bar{u}$ of its members and form pointwise $\gamma$-credible intervals $C_{\gamma} := [q_{(1-\gamma)/2},\, q_{(1+\gamma)/2}]$, where $q_p(\boldsymbol{x})$ denotes the empirical $p$-quantile of the member values $\{ u_{(s)}(\boldsymbol{x}) \}_{s=1}^{S}$.

\paragraph{Relative $\ell^1$-errors}

For a nodal vector $\mathbf{v}$, we define the relative $\ell^1$-error as $E(\mathbf{v}) := \| \mathbf{v} - \mathbf{u}^{\dagger} \|_1 / \| \mathbf{u}^{\dagger} \|_1$ and report two statistics:
(i) the relative $\ell^1$-error of the ensemble mean, $E(\bar{\mathbf{u}})$;
and (ii) the mean of the relative $\ell^1$-errors, $\frac{1}{S} \sum_{s=1}^{S} E(\mathbf{u}_{(s)})$.
The first statistic measures the accuracy of the point estimate, whereas the second characterizes the accuracy of the individual ensemble members.
By the triangle inequality, the second is never smaller than the first, and their gap indicates how much the ensemble benefits from the diversity of its members.

\paragraph{Coverage test}

To validate the uncertainty estimates, we consider the frequentist coverage \cite{lakshminarayanan2017simple,hoffmann2021deep} of the credible intervals $C_{\gamma}$.
The empirical coverage at nominal level $\gamma$ is $\mathrm{cov}(\gamma) := \frac{1}{N} \sum_{n=1}^{N} \mathbbm{1}\{ u^{\dagger}(\mathbf{x}_n) \in C_{\gamma}(\mathbf{x}_n) \}$.
Plotting $\mathrm{cov}(\gamma)$ against the nominal levels $\gamma \in [0.5, 0.99]$ yields a calibration curve, where curves on the diagonal indicate well-calibrated uncertainty, while curves below (above) the diagonal indicate overconfident (conservative) intervals.

\paragraph{Continuous ranked probability score}

We further evaluate the deep ensembles' uncertainty estimates using the continuous ranked probability score (CRPS).
The CRPS at position $\boldsymbol{x}$ is $\text{crps}(\boldsymbol{x}) = \frac{1}{S} \sum_{s=1}^{S} | u_{(s)}(\boldsymbol{x}) - u^\dagger(\boldsymbol{x}) | - \frac{1}{2S^2} \sum_{s,s'=1}^{S} | u_{(s)}(\boldsymbol{x}) - u_{(s')}(\boldsymbol{x}) |$.
Averaging over the domain yields $\text{CRPS} = \frac{1}{|\Omega|} \int_\Omega \text{crps}(\boldsymbol{x}) \,\mathrm{d}\boldsymbol{x} \approx \frac{1}{N} \sum_{n=1}^{N} \text{crps}(\mathbf{x}_n)$.
The CRPS generalizes the absolute error, to which it reduces for $S = 1$. 
The first term measures how far the members lie from the truth, while the second rewards their spread.
It is thus small only where the members concentrate near $u^{\dagger}$. 
Where the ensemble misses $u^{\dagger}$, it penalizes overconfident ensembles.

\subsection{Denoising a piecewise constant signal}
\label{subsec:num_denoising}

\begin{figure}
    \centering
    \includegraphics[width=0.99\textwidth]{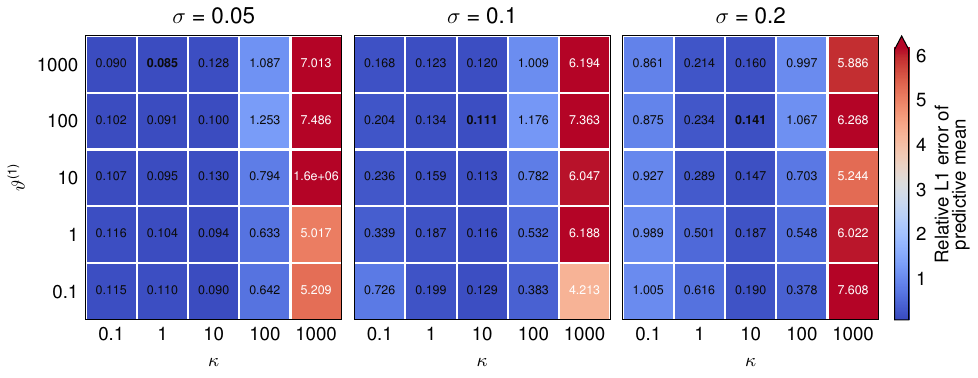}
	\caption{
		The influence of the $\kappa$ and $\vartheta^{(1)}$ parameters for different noise levels.
		We find that $\kappa=10$ and $\vartheta^{(1)}=100$ form a reliably well-performing parameter pair across all noise levels.
	}
	\label{fig:Denoising_ParameterChoice}
\end{figure}

We consider recovering a piecewise constant signal $u: \Omega = [-1, 1] \to \R$ based on $M=50$ noisy direct observations with additive observational noise $\mathbf{e} \sim \mathcal{N}(0, \sigma^2)$ and $\sigma = 0.1$.
Unless stated otherwise, we use the GSBL-BNN prior with $r=-1$, $\beta=1$, $\vartheta^{(1)}=100$, $c=2$, and $\kappa=10$.
We choose $r=-1$ and $\beta=1$ because these are typical choices in the GSBL literature \cite{glaubitz2023generalized,xiao2023sequential,glaubitz2024leveraging,lindbloom2025efficient}.
While we made no attempt to optimize the parameters of GSBL-BNN prior, we found their performance to be relatively robust w.\,r.\,t.\ them. 
We later demonstrate this for $\vartheta^{(1)}$ and $\kappa=10$ in \cref{fig:Denoising_ParameterChoice} but omit further sensitivity studies for brevity. 
If not stated otherwise, we consider the GSBL-BNN with shared-$\boldsymbol{\psi}$ biases (see \Cref{subsubsec:prior_spec_bias}).

We start by comparing our GSBL-BNNs with Gaussian- and Cauchy-BNNs that respectively put a standard normal prior and a centered Cauchy prior with scale $0.1$ on both the weight and bias parameters \cite{li2022bayesian}. 
While this comparison is not exhaustive, it highlights some distinct differences between our GSBL-BNNs and the other two. 
In all cases, we use a neural network with four hidden layers, each with $64$ neurons.
Furthermore, we train with the Muon optimizer for $500$ iterations and then refine with L-BFGS for $1500$ iterations.

\begin{figure}[tb]
	\centering
	\includegraphics[width=\textwidth]{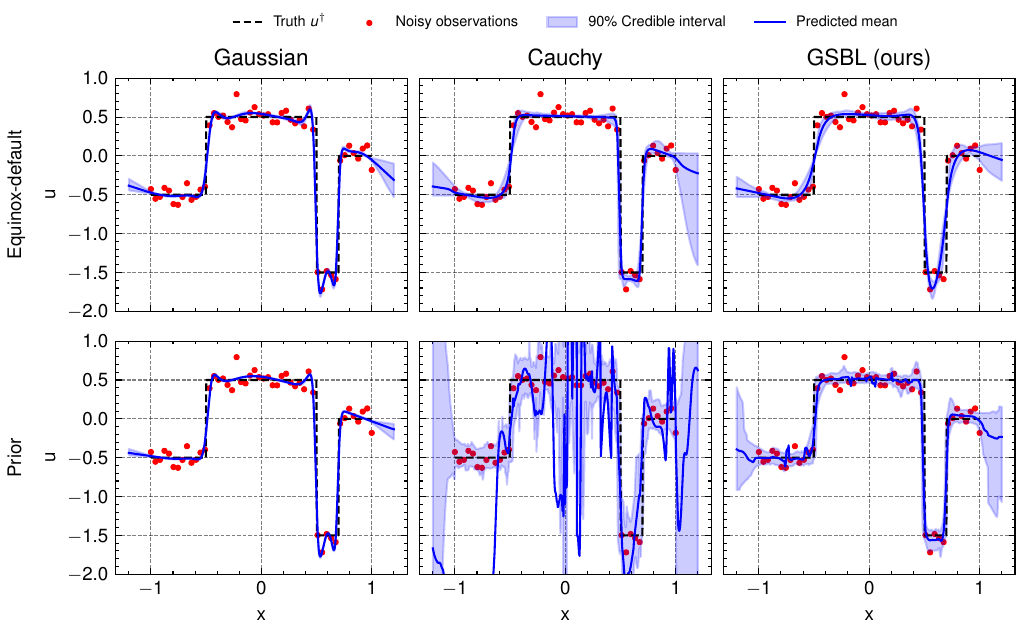}
	\caption{
		Denoising a piecewise constant signal.
		We compare deep ensembles ($S=50$) for a Gaussian (left), Cauchy (center), and our GSBL-BNN (right, with shared-$\boldsymbol{\psi}$).
		We randomly initialize the BNN ensembles using the Equinox-default (top) and draws from the respective prior (bottom).
	}\label{fig:Denoising_CompareBNNs}
\end{figure}

\Cref{fig:Denoising_CompareBNNs} illustrates the piecewise constant signal, the observational data, and the results (predicted mean and $90\%$ credible intervals) of a deep ensemble of size $S=50$ for all three BNNs.
Furthermore, we compare randomly initializing each BNN using the Equinox default\footnote{%
The Equinox-default initialization \cite{kidger2021equinox} draws all weights and biases of layer $\ell$ i.i.d.\ from $\mathcal{U}(-n_{\ell-1}^{-1/2}, n_{\ell-1}^{-1/2})$, i.e., with variance $1/(3 n_{\ell-1})$, a factor of three below LeCun scaling~\cite{lecun_efficient_1998}.
} (top) and drawing from the Gaussian, Cauchy, and GSBL priors, respectively. 
The Gaussian-BNN shows oscillatory behavior around the jump discontinuities, demonstrating its failure to capture large gradients.
At the same time, we observe relatively small credible intervals for the Gaussian-BNNs, even in regions where they are visibly inaccurate, making them confidently wrong.
The Cauchy-BNN with Equinox-default initialization accurately recovers the sharp gradient, but shows almost no difference between different members.
When initializing the Cauchy-BNN with random prior draws, however, the ensemble mean becomes highly oscillatory and inaccurate.
The GSBL-BNN with Equinox-default initialization yields a reconstruction comparable to the Cauchy-BNN, although it captures the constant region $0.5 \leq x \leq 0.7$ less accurately.
When using the prior initialization, the GSBL-BNN accurately predicts the sharp gradients.
Overall, we observe the GSBL-BNN to be more robust w.\,r.\,t.\ the initialization than the Cauchy-BNN.
From here on, we consider only GSBL-BNNs with prior initialization and Gaussian- and Cauchy-BNNs with Equinox default initialization.

\begin{figure}[tb]
	\centering
	\includegraphics[width=\textwidth]{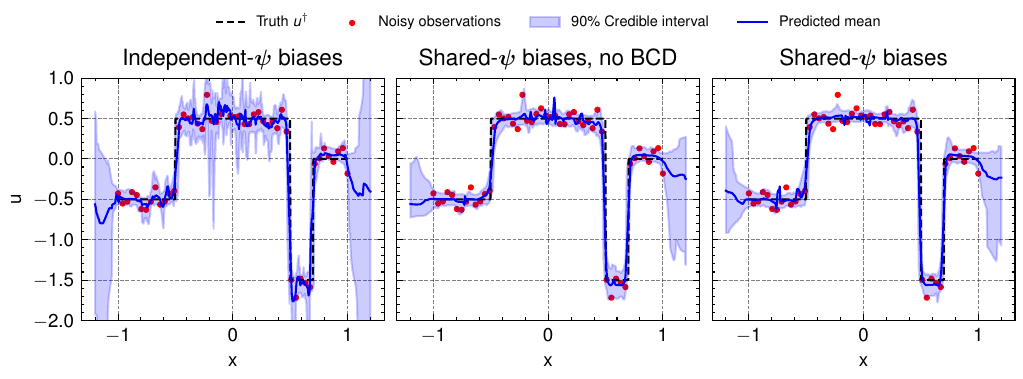}
	\caption{
		Denoising a piecewise constant signal. 
		GSBL-BNNs with independent-$\boldsymbol{\psi}$ biases and BCD (left), shared-$\boldsymbol{\psi}$ biases and no BCD (center), and shared-$\boldsymbol{\psi}$ biases and BCD (right).
	}
	\label{fig:Denoising_GSBLComparison}
\end{figure}

We next study the influence of our shared-$\psi$ bias contribution and the BCD optimization on the GSBL-BNNs.
\Cref{fig:Denoising_GSBLComparison} compares GSBL-BNNs with independent-$\boldsymbol{\psi}$ biases and BCD (left), shared-$\boldsymbol{\psi}$ biases and no BCD (center), and shared-$\boldsymbol{\psi}$ biases and BCD (right).
No BCD means that the neural network parameter $\boldsymbol{\theta}$ and the hyperparameter $\boldsymbol{\psi}$ are trained jointly---here using Muon and L-BFGS optimization.
The comparison shows that the shared-$\boldsymbol{\psi}$ biases (joint and BCD training) make the GSBL-BNN approximations less oscillatory overall.
The BCD optimization procedure reduces oscillations even further, making it the preferred choice.

\cref{tab:denoising_errors} lists the relative $\ell^1$-error of the predictive mean, the CRPS, and the mean relative $\ell^1$-error of the ensemble members for the different BNN priors.
The predictive mean reflects the overall accuracy of the deep ensemble and points to the Cauchy-BNN and GSBL-BNN with shared-$\boldsymbol{\psi}$ as the best choices.
However, CRPS is the most meaningful value for our purposes because it balances accuracy with confidence.
Considering this, our GSBL-BNN is the best choice.
In terms of the mean relative $\ell^1$-error of the ensemble members, the Cauchy-BNN appears stronger than our GSBL-BNN.
However, this error differs only slightly from that of its predictive mean ($1.28\mathrm{e}{-1}$ vs.\ $1.16\mathrm{e}{-1}$).
One potential reason is that all members of the Cauchy-BNN ensemble collapse into the (almost) same local minimum.

\begin{table}[tb]
	\centering
	\caption{
		Relative $\ell^1$-error of the predictive mean (predictive), CRPS, and mean relative $\ell^1$-error of the ensemble members (mean) for the denoising problem.
		Best values are in bold.
	}
	\label{tab:denoising_errors}
	\small
	\begin{tabular}{ll*{3}{c}}
		\toprule
		prior & initialization & predictive & CRPS & mean \\
		\midrule
		\multirow{2}{*}{Gaussian}
			& Equinox-default & $1.28\mathrm{e}{-1}$ & $6.13\mathrm{e}{-2}$ & $1.29\mathrm{e}{-1}$ \\
			& prior draw    & $1.24\mathrm{e}{-1}$ & $6.25\mathrm{e}{-2}$ & $\mathbf{1.25\mathrm{\textbf{e}}{-1}}$ \\
		\addlinespace[2pt]
		\multirow{2}{*}{Cauchy}
			& Equinox-default & $1.16\mathrm{e}{-1}$ & $4.77\mathrm{e}{-2}$ & $1.28\mathrm{e}{-1}$ \\
			& prior draw    & $1.49\mathrm{e}{+0}$ & $7.75\mathrm{e}{-2}$ & $1.93\mathrm{e}{+0}$ \\
		\addlinespace[2pt]
		\multirow{2}{*}{GSBL (ours)}
			& Equinox-default & $1.92\mathrm{e}{-1}$ & $8.42\mathrm{e}{-2}$ & $2.08\mathrm{e}{-1}$ \\
			& prior draw    & $\mathbf{1.11\mathrm{\textbf{e}}{-1}}$ & $\mathbf{3.78\mathrm{\textbf{e}}{-2}}$ & $1.66\mathrm{e}{-1}$ \\
		\bottomrule
	\end{tabular}
\end{table}

\begin{figure}
	\centering
	\begin{subfigure}[b]{.48\textwidth}
		\centering
		\includegraphics[height=.716\linewidth]{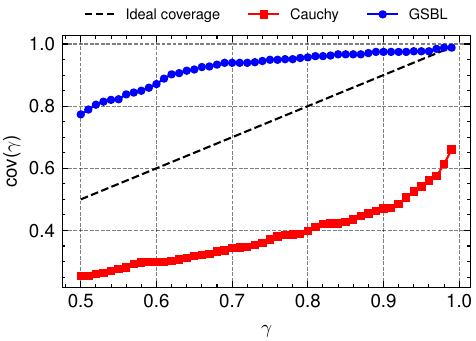}
		\caption{Coverage test}\label{fig:Denoising_Coverage}
	\end{subfigure}
	~
	\begin{subfigure}[b]{.48\textwidth}
		\centering
		\includegraphics[height=.716\linewidth]{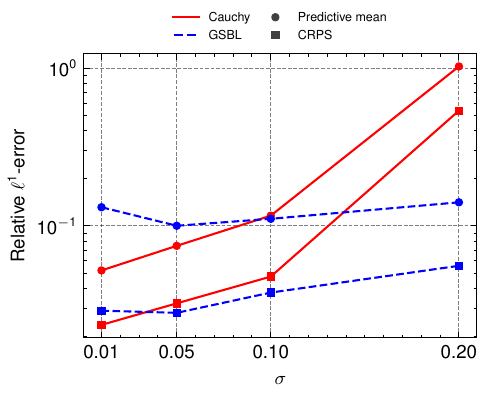}
		\caption{Accuracy versus noise level}\label{fig:Denoising_NoiseComparison}
	\end{subfigure}
	\caption{
		Comparison of the Cauchy-BNN and our GSBL-BNN for the signal denoising problem.
		(a) Coverage test: the Cauchy-BNN underestimates uncertainty, while our GSBL-BNN overestimates it but is closer to the ideal coverage.
		(b) Relative $\ell^1$-error of the predictive mean (circles) and CRPS (squares) for the Cauchy-BNN (solid) and our GSBL-BNN (dashed) over the noise level $\sigma$.
	}\label{fig:Denoising_CauchyGSBL}
\end{figure}

To better understand the estimated uncertainty of different BNNs, \cref{fig:Denoising_Coverage} illustrates the coverage for the Cauchy-BNN and our GSBL-BNN.
While the Cauchy-BNN underestimates uncertainty, our GSBL-BNN overestimates it.
At the same time, we observe that our GSBL-BNN provides better calibrated uncertainty estimates than the Cauchy-BNN in the sense that its empirical coverage is closer to the nominal levels, especially for larger $\gamma$. 
Specifically, the mean-squared deviation between empirical and ideal coverage, $\mathrm{MSE}_{\mathrm{cov}} = \frac{1}{I}\sum_{i=1}^I ( \mathrm{cov}(\gamma_i) - \gamma_i )^2$, is $4.05\mathrm{e}{-2}$ for our GSBL-BNN and $1.34\mathrm{e}{-1}$ for the Cauchy-BNN.

Moreover, \Cref{fig:Denoising_NoiseComparison} illustrates the behavior of our GSBL-BNN and the Cauchy-BNN across different noise levels.
While we previously considered $\sigma = 0.1$, 
\Cref{fig:Denoising_NoiseComparison} shows that our GSBL-BNN is notably more accurate than the Cauchy-BNN at higher noise levels ($\sigma \geq 0.1$), both in terms of the relative $\ell^1$-error of the predictive mean and the CRPS. 
The Cauchy-BNN becomes more competitive at low noise levels ($\sigma \leq 0.05$).

The GSBL-BNN prior has several free parameters, including $\vartheta^{(1)}$ and $\kappa$.
\cref{fig:Denoising_ParameterChoice} thus provides a brief parameter study for $r=-1.0,\ \beta=1.0$, and $c=2$.
We report the relative $\ell^1$-error of the predictive mean for the GSBL-BNN for different noise levels $\sigma \in \{0.05, 0.1, 0.2\}$ and parameters $\vartheta^{(1)}, \kappa \in \{0.1, 1, 10, 100, 1000\}$.
We observe that the combination of $\kappa = 10$ and $\vartheta^{(1)} = 100$ performs well across all considered noise levels.
Furthermore, the error of our GSBL-BNN is robust for large ranges of $\kappa$ and $\vartheta^{(1)}$.

\subsection{Deblurring a piecewise smooth signal}
\label{subsec:num_deblurring}

Consider recovering a piecewise smooth signal $u:[0,1] \to \R$ with jump discontinuities from noisy blurry observations. 
The corresponding forward operator is $F_{i,n} = \exp( - [\tilde{\mathbf{x}}_i - \mathbf{x}_n)^2] / 2 \omega^2 ) / \sum_{k=1}^{N} \exp( - [\tilde{\mathbf{x}}_i - \mathbf{x}_n)^2] / 2 \omega^2 )$, where $\tilde{\mathbf{x}}_i$ is the $i$th observation point, $\mathbf{x}_n$ is the $n$th of the $N$ grid points, and $\omega$ denotes the blur width, which we set to $3 \cdot 10^{-2}$.
We discretize the computational domain $\Omega = [0, 1]$ using $N=128$ equidistant points and observe the blurred signal $F\mathbf{u}^\dagger$ with additive observational noise $\mathbf{e} \sim \mathcal{N}(0, \sigma^2)$ and $\sigma = 0.1$ at $M=50$ points.
We again use $r=-1$, $\beta=1$, $\vartheta^{(1)}=100$, $c=2$, and $\kappa=10$ for the prior parameters of the GSBL-BNN.
The initial training with the Muon optimizer runs for $500$ iterations, and the L-BFGS refinement runs for $1500$ iterations.
As before, we have made no attempt to optimize the GSBL-BNN parameters, the architecture, or the training for this problem.
We again compare our GSBL-BNN with BNNs using a standard normal prior and a centered Cauchy prior with scale $0.1$.
GSBL-BNNs are always initialized using random prior draws, while Gaussian-BNNs and Cauchy-BNNs are initialized using the Equinox default.

\begin{figure}[tb]
	\centering
	\includegraphics[width=\textwidth]{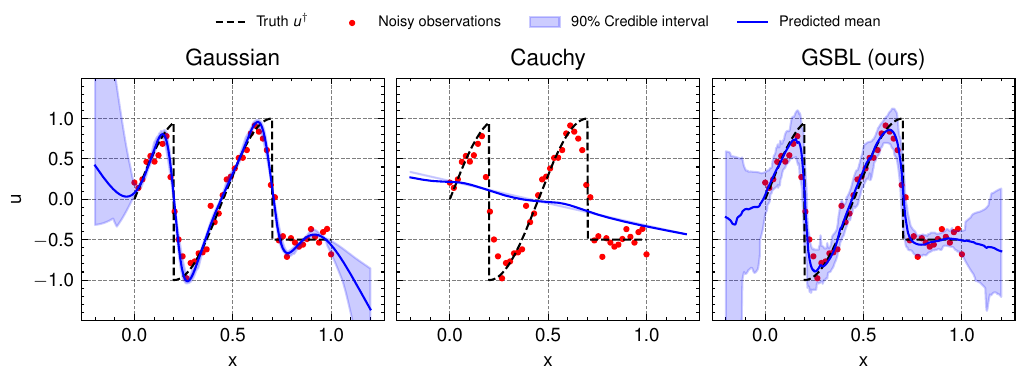}
	\\[-5pt]
	\caption{
		Deblurring a piecewise smooth signal.
		We compare deep ensembles (size $S=50$) for a Gaussian (left), a Cauchy (center), and our GSBL-BNN (right, with shared-$\boldsymbol{\psi}$).
	}
	\label{fig:Deblurring_CompareBNNs}
\end{figure}

\cref{fig:Deblurring_CompareBNNs} illustrates the piecewise smooth signal, the observational data, and the predictive mean and 90\% credible intervals of the Gaussian-BNN, Cauchy-BNN, and our GSBL-BNN.
The Gaussian-BNN blurs the edges, does not capture the constant region for $0.7 < x \leq 1.0$, and again yields a collapsed ensemble.
The Cauchy-BNN produces no meaningful solution at all.
Our GSBL-BNN also blurs the edges, but better captures the constant region, and yields more meaningful credible intervals that engulf the ground truth.

\cref{tab:deblurring_errors} lists the relative $\ell^1$-error of the predictive mean, the CRPS, and the mean relative $\ell^1$-error of the ensemble members for the different BNN priors. 
As in \Cref{subsec:num_denoising}, our GSBL-BNN yields the smallest relative $\ell^1$-error of the predictive mean and CRPS.

\begin{table}[tb]
	\centering
	\caption{
		Relative $\ell^1$-error of the predictive mean (predictive), CRPS, and mean relative $\ell^1$-error of the ensemble members (mean) for the deblurring problem.
		Best values are in bold.
	}\label{tab:deblurring_errors}
	\small
	\begin{tabular}{l*{3}{c}}
		\toprule
		prior & predictive & CRPS & mean \\
		\midrule
		Gaussian & $1.97\mathrm{e}{-1}$ & $1.04\mathrm{e}{-1}$ & $\mathbf{2.08\mathrm{\textbf{e}}{-1}}$ \\
		Cauchy   & $8.63\mathrm{e}{-1}$ & $4.88\mathrm{e}{-1}$ & $8.72\mathrm{e}{-1}$ \\
		GSBL (ours) & $\mathbf{1.48\mathrm{\textbf{e}}{-1}}$ & $\mathbf{6.90\mathrm{\textbf{e}}{-2}}$ & $2.49\mathrm{e}{-1}$ \\
		\bottomrule
	\end{tabular}
\end{table}

\subsection{Inverting a two-dimensional hyperbolic conservation law}
\label{subsec:num_2d}

To demonstrate that the proposed GSBL-BNNs extend in a straightforward manner to more than one dimension and non-linear problems, we consider the two-dimensional Kurganov--Petrova--Popov (KPP) rotating wave problem \cite{kurganov2007adaptive}: 
\begin{equation}
	u_t + \left(\sin u\right)_x + \left(\cos u\right)_y = 0, 
	\quad 
	u(x, y, 0) = \begin{cases}
		3.5 \pi, & \text{if } x^2 + y^2 < 1, \\
		0.25 \pi, & \text{otherwise},
	\end{cases}
\end{equation}
on $(x, y) \in \Omega = [-2, 2]\times[-2.5, 1.5]$ with transmissive boundary conditions. 
Our goal is to reconstruct the piecewise constant, disk-shaped initial state $u(x, y, 0)$ from noisy sparse observations of the solution at $t=1$ on a $22 \times 22$ lattice of equally spaced points ($M = 484$ observations).
The forward operator is the nonlinear numerical solution operator, mapping an initial state to the solution at $t=1$.
We discretize the domain $\Omega$ into $64 \times 64$ equally sized rectangular elements and solve the PDE numerically with a discontinuous Galerkin method \cite{hesthaven2007nodal,cockburn2012discontinuous}.
To this end, we use polynomials of degree one, a local Lax--Friedrichs numerical flux, a minmod slope limiter, and a third-order three-stage strong stability preserving Runge--Kutta method with CFL number $0.9$. 
To avoid the ``inverse crime''~\cite{kaipio2007statistical}, we generate the sparse observations at $t=1$ using a refined grid of $512 \times 512$ elements and add zero-mean Gaussian noise with standard deviation $\sigma = 0.5$.

We use a GSBL-BNN with a two-dimensional input, four hidden layers with $128$ neurons each, and a scalar output.
We again use $r=-1,\ \beta=1,\ c=2$, and $\kappa = 10$.
Furthermore, we set $\vartheta^{(1)} = 10$, as we found this value to perform better than $\vartheta^{(1)} = 100$, which we used in the previous two examples. 
Notably, the range of the unknown function value was $[-1.5, 0.5]$ and $[-1, 1]$ in the previous examples, but here is $[0.25\pi, 3.5\pi]$.
To account for the changed range, we apply an affine transformation to the output of our GSBL-BNN that maps $[0, 1]$ onto $[0.25\pi, 3.5\pi]$.
We compute a deep ensemble with $S=10$ members.
We first train for $2500$ steps with the Muon optimizer at a fixed learning rate of $10^{-3}$, followed by $12500$ refinement steps with L-BFGS.
One of the ten runs terminates at a training loss roughly seven times higher than the remaining nine and is discarded as a failed optimization.  
For brevity, we refrain from comparing the GSBL-BNN to existing Gaussian and Cauchy BNNs.

\begin{figure}[tb]
	\centering
	\includegraphics[width=\textwidth]{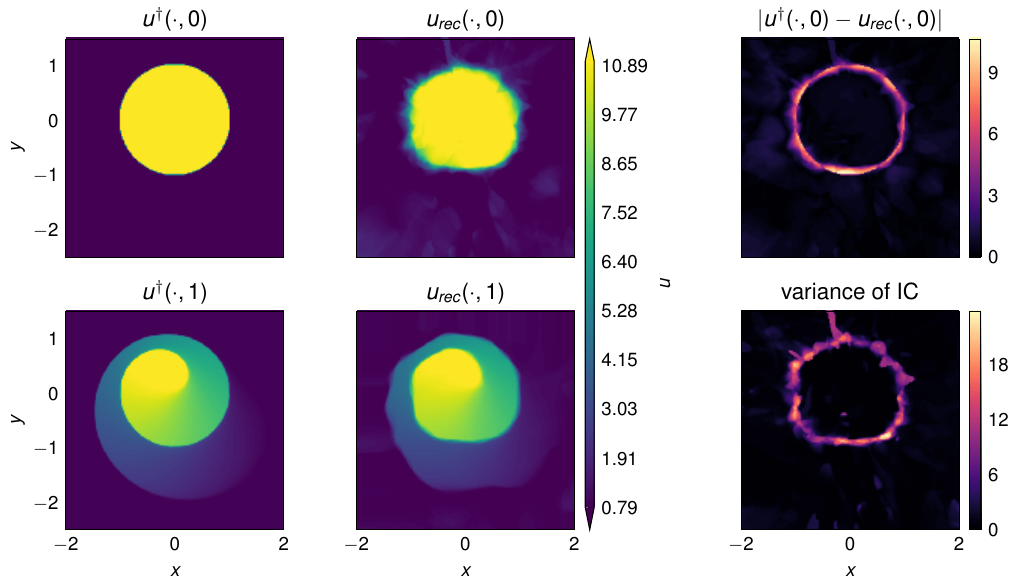}
	\caption{
		GSBL-BNN reconstruction of the piecewise constant initial state of the KPP rotating wave problem.
		Top: The true initial state (top, left), the predictive mean of the recovered initial state (top, center), and their absolute difference (top, right). 
		Bottom: The true solution at $t=1$ (bottom, left), the solution at $t=1$ computed from the recovered initial state (bottom, center), and the variance of the recovered initial state (bottom, right).
	}\label{fig:2D_EnsemblePlots}
\end{figure}

\cref{fig:2D_EnsemblePlots} illustrates the true and recovered initial states, the corresponding solutions at $t=1$, the absolute error of the recovered initial state, and its ensemble variance.
The initial state recovered by the GSBL-BNN matches the ground truth relatively closely, given the difficulty of the inverse problem.
This is also reflected in the relative $\ell^1$-error of the predictive mean of the initial state, which is $E(\bar{\mathbf{u}}) = 2.11\mathrm{e}{-1}$. 
Evolving this initial state to $t=1$ results in a recovered solution with relative $\ell^1$-error of $E(\mathbf{\bar{u}(\cdot, 1)}) = 7.80\mathrm{e}{-2}$.
\Cref{fig:2D_EnsemblePlots} further illustrates that the errors at the edges dominate the overall error.
We also observe the variance to be concentrated along the edges, illustrating that our GSBL-BNN is more uncertain precisely where it is wrong---thus avoiding confidently wrong predictions.
The faint streak of non-zero variance above the disk originates from a single spurious prediction.

\begin{figure}[tb]
	\centering
	\includegraphics[width=\textwidth]{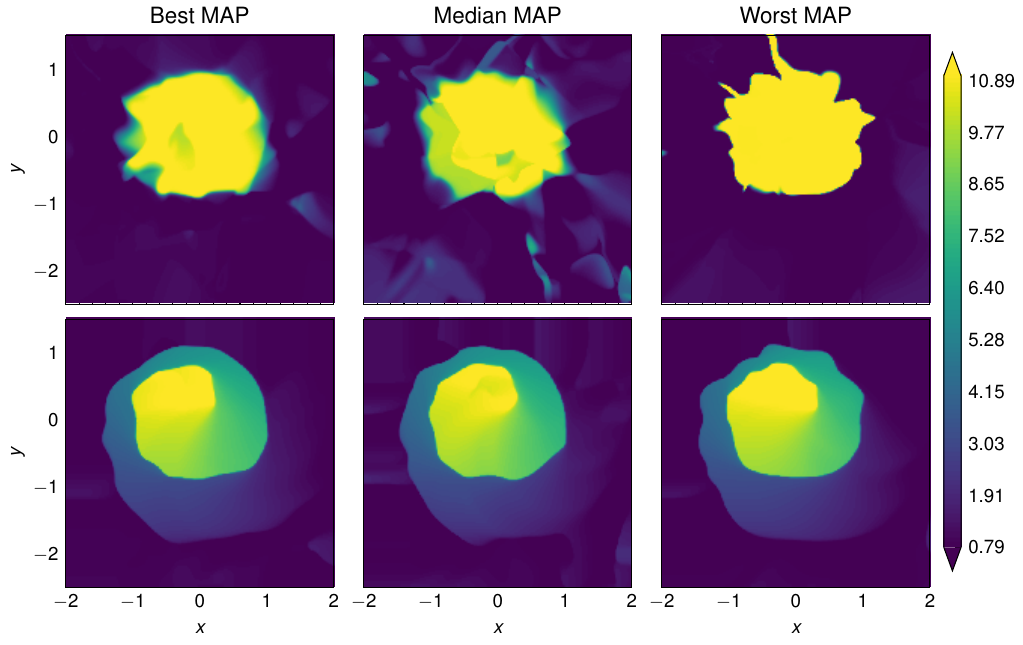}
	\caption{
		The GSBL-BNN ensemble member with the lowest (left), median (center), and highest (right) training loss for the KPP problem.
		Initial states (top) and their solutions at $t=1$ (bottom).
	}\label{fig:2D_EnsembleMember}
\end{figure}

To demonstrate the variability among the ensemble members, \cref{fig:2D_EnsembleMember} illustrates the ensemble members with the lowest, median, and highest training loss.
Although their numerical solutions at $t=1$ are relatively similar, the recovered initial states differ considerably---highlighting the ill-posedness of the considered inverse problem.

\section{Concluding thoughts}
\label{sec:summary}

We introduced GSBL-BNNs, which model the neural parameters by conditionally Gaussian priors with generalized gamma hyperpriors and include Laplace, Student's $t$, and Cauchy marginal priors as special cases.
Unlike existing edge-preserving BNNs, GSBL-BNNs couple each neuron's bias to its weights through a shared hyperparameter, which spreads the induced jumps across the domain while keeping the biases conditionally Gaussian.
Conditional Gaussianity enables our BCD approach for MAP estimation, which alternates between network training with iteratively reweighted $\ell^2$-regularization and exact hyperparameter updates.
Deep ensembles of these MAP estimates provide approximate UQ.
In denoising and deblurring experiments, GSBL-BNNs yielded the most accurate predictive means and the lowest CRPS among the compared BNNs.
They were also more robust than Cauchy-BNNs to initialization and noise.
For a two-dimensional hyperbolic conservation law, the ensemble variance concentrated along the edges, where the errors were largest.

Several questions remain open.
Deep ensembles are a heuristic, and we will pursue more principled MCMC or variational inference in future work.
Partial BNNs \cite{daxberger2021bayesian,sharma2023bayesian,zhao2024feynman}, which treat only some neural parameters as random, could reduce the computational cost.
Finally, GSBL-BNNs could provide accurate, uncertainty-aware representations of discontinuous states in neural PDE solvers \cite{mao2020physics,yang2021bpinns} and data assimilation for hyperbolic conservation laws \cite{zhou2026neural,chandravamsi2026feature,glaubitz2026regularity}.

\section*{Acknowledgements}
JG acknowledges support from the Swedish Research Council (VR) Starting Grant \#2025-05370, the Zenith Career Development Grant \#26.07, the Swedish e-Science Research Centre (SeRC), and the National Academic Infrastructure for Supercomputing in Sweden (NAISS) grant \#2025/22-1599.
\section*{Declaration}

Generative AI (Gemini and Claude Code) was used to assist with code and writing. 
All AI-generated output was tested and reviewed by the authors.


\bibliographystyle{siamplain}
\bibliography{references}

\end{document}